\documentclass[11pt]{article}
\usepackage{amsfonts}
\usepackage{amsthm}
\usepackage{amsmath}
\usepackage{amssymb}
\usepackage{amscd}
\usepackage{verbatim}
\usepackage{multicol}
\usepackage{tabularx}

\theoremstyle{definition}
\newtheorem{thm}{Theorem}[section]
\newtheorem{prop}[thm]{Proposition}
\newtheorem{cor}[thm]{Corollary}
\newtheorem{lem}[thm]{Lemma}
\newtheorem{defn}{Definition}[section]

\newtheorem{Thm}{Theorem}

\numberwithin{equation}{section}

\def\obs{\noindent{\normalsize{\bf Remark: }}}

\def\blb#1{\text{$\mathbb{#1}$}}
\def\cal#1{\text{$\mathcal{#1}$}}

\def\ord#1^#2{#1$^{\text{#2}}$}

\def\lie#1{\mathfrak{#1}}

\def\hlie#1{\widehat{\mathfrak{#1}}}

\def\dnl{=\hspace{-.2cm}>\hspace{-.2cm}=}
\def\dnr{=\hspace{-.2cm}<\hspace{-.2cm}=}
\def\sn{-\hspace{-.17cm}-}

\begin{document}

\flushbottom

\author{Pavel I. Etingof\footnote{MIT - 77 Massachussets Av., Cambridge MA 02139 USA - 
Room 2-176} \hspace{1.15cm} Adriano A. Moura\footnote{IMECC/UNICAMP - Caixa Postal: 6065, CEP: 13083-970 - 
Campinas SP Brazil}\\   {\small etingof@math.mit.edu \hspace{.5cm} adrianoam@ime.unicamp.br }}
\title{Elliptic Central Characters and Blocks \\ of Finite
  Dimensional Representations\\ of Quantum Affine Algebras}
\date{April/2002}
\maketitle

\nocite{AkKa97}
\nocite{Bec93}
\nocite{Cha01}
\nocite{ChPr91,ChPr93,ChPr94bo,ChPr95,ChPr96}
\nocite{Dri88}
\nocite{EFK98}
\nocite{EtMo01}
\nocite{FreMu01}
\nocite{FreRe98,FreRe00}
\nocite{Kas00}
\nocite{KaSo95}
\nocite{KhTo92}
\nocite{Mou01}


\centerline{\bf For Igor Frenkel, on the occasion of his 50-th birthday}

\section*{Introduction}

In this paper we describe the block decomposition
of the category of finite dimensional representations 
of a quantum affine algebra $U_q(\hlie g)$, where $|q|< 1$.
Namely, we find that the blocks are parametrized by 
elliptic central characters, which are certain elliptic
functions attached to irreducible representations. 

The plan of the paper is as follows. 

In Section 1, we recall the basics about blocks in abelian
categories. 

In Section 2, we recall the definition of $U_q(\hlie g)$ and
the basic facts about its finite dimensional representations. 

In Section 3, we define the elliptic central character
of finite dimensional representations of $U_q(\hlie g)$.
Namely, we show that if $X,Y$ are such representations 
and $X$ is irreducible, then the operator 
$R^{2,1}_{Y,X}(z^{-1})R_{X,Y}(z): X\otimes Y\to X\otimes Y$
(where $R$ is the R-matrix) is of the form $1\otimes
\xi_X(z)|_Y$, where $\xi_X(z)$ is an endomorphism of the identity functor 
with coefficients in the field of elliptic functions. 
If $Y$ is also irreducible, then $\xi_X(z)|_Y$ 
is a scalar. The collection of these scalars 
for all irreducible $X$ is called the elliptic central character of
$Y$. (Thus, the elliptic central character plays the role of 
non-existent nontrivial central elements of $U_q(\lie g)$). Our main result is

\begin{Thm}\label{mainth} Blocks in the category of finite dimensional
representations of $U_q(\hlie g)$ consist of representations
whose simple constituents have a given elliptic central
character.
\end{Thm}

The rest of the paper is devoted to the proof of Theorem \ref{mainth}. 
This proof is based on a result of Chari and Kashiwara
on the cyclicity of tensor products of fundamental
representations, and a lengthy case-by-case analysis. 
It would be interesting to obtain a uniform proof. 

In Section 4, we prove Theorem \ref{mainth} for type A.

In Section 5, we recall the Drinfeld realization of
$U_q(\hlie g)$, 
and explain how to compute elliptic central characters of
fundamental representations. 

In Section 6, we prove Theorem \ref{mainth} for types B-G. This is 
much more technically challenging than type A
(for types $E$ and $F$ we relied on a 
computer to perform some calculations).   

In the appendix, we list the formulas used in Section 6
to compute the elliptic central characters; 
some of them were obtained using a computer. 

We note that although elliptic central characters were used in
this paper for a particular purpose (to classify blocks of finite
dimensional representations), they may be used in other
problems about $U_q(\hlie g)$. For example, as is shown below, 
they can sometimes be used to decide when an irreducible 
finite dimensional representation
occurs in the tensor product of two others. 
Therefore, we feel that elliptic central characters are worthy of
further study. In particular, it would be interesting to study 
their connections with other objects in representation theory of
$U_q(\hlie g)$, such as minimal affinizations and q-characters.

\noindent{\bf Acknowledgments: }
This paper is dedicated to the 50-th birthday of Igor Frenkel. 
The first author thanks Igor Frenkel for being a wonderful
adviser, and in particular for introducing him to the circle of
ideas that led to this paper and to many others. 

We are grateful to Edward Frenkel and Vyjayanthi Chari  for
very useful discussions and explanations, which played a crucial role in the creation of
 this paper 
(in particular, E. Frenkel explained to us how to compute the function $f_{X,Y}(z)$ of proposition \ref{p:factor}). 
We also thank E. Mukhin for a useful discussion. The work of P.E. was partially supported by the NSF
grant DMS-9988796, and 
was done in part for the Clay Mathematics Institute. A.M. is grateful to MIT for hospitality. The 
Ph.D. studies of A.M., during which this work was accomplished,  
are supported by FAPESP (99/11600-0), and his visit to MIT
is supported by CAPES (0365/01-7), Brazil. 


\section{Block decomposition of an abelian category}

Let us recall the basics about blocks in abelian categories. 
This material is standard, and we give it for the reader's convenience.

Let \cal C be an abelian category, in which every object has
finite length. In this case, it is well known 
that any object is uniquely representable 
as a direct sum of indecomposable objects.  

\begin{defn}
Two indecomposable objects
$X_1,X_2$ of \cal C are {\it linked} if there is no splitting 
of \cal C in a direct sum of two abelian categories, 
$\mathcal C=\mathcal C_1\oplus \mathcal C_2$, such that 
$X_1\in \mathcal C_1$ and $X_2\in \mathcal C_2$. 
\end{defn}

It is easy to see that linking is an equivalence relation. 

\begin{prop}\label{blockdec} The category $\mathcal C$ admits a unique
decomposition into a direct sum of indecomposable abelain
categories: $\mathcal C=\oplus_{\alpha\in I}\mathcal
C_\alpha$.  
\end{prop}

\begin{proof} Let $I$ be the set of equivalence classes of linked
indecomposable objects, and for $\alpha\in I$ let
$\mathcal C_\alpha$ be the subcategory of $\mathcal C$, 
consisting of direct sums of objects from $\alpha$. 
By the uniqueness of the decomposition into indecomposables, 
we have $\mathcal C=\oplus_{\alpha\in I} \mathcal C_\alpha$. 
Furthermore, the categories $\mathcal C_\alpha$
are indecomposable. Indeed, if $\mathcal C_\alpha=\mathcal
C_\alpha^1\oplus \mathcal C_\alpha^2$ is a nontrivial
decomposition, then 
any indecomposables $X_1\in \mathcal C_\alpha^1$, 
 $X_2\in \mathcal C_\alpha^2$ are not linked -- contradiction. 

The uniqueness of the decomposition is obvious.
\end{proof}

\begin{defn} The subcategories $\mathcal C_\alpha$ are called the
  {\it blocks} of $\mathcal C$, and the decomposition of 
Proposition \ref{blockdec} is called the {\it block
  decomposition} of $\mathcal C$. 
\end{defn}

Recall that for any $X\in \mathcal C$, one can uniquely 
specify the simple objects (with multiplicities) which occur
as constituents in $X$ (the Jordan-Holder Theorem).  

The following trivial lemma will be used below. 

\begin{lem}\label{l:linkcond}\hfill\vspace{-.25cm}
\begin{enumerate}
\item Two simple objects are linked if they occur as constituents
of the same indecomposable object.
\item Two indecomposable objects are linked if they have some linked simple constituents.
\end{enumerate}
\end{lem}

\begin{proof}
Clear.
\end{proof}

\section{Quantum Affine Algebras and the Category \cal C}

Let $\lie g$ be a finite-dimensional complex simple Lie algebra
with Cartan subalgebra $\lie h$. The simple roots will be denoted
$\alpha_1,\dots, \alpha_n$, the fundamental weights
$\omega_1,\dots, \omega_n$ and the invariant bilinear form $<,>$
is normalized so
that $<\theta,\theta>=2$ for the maximal root $\theta=\sum \theta_i\alpha_i$.
The Cartan matrix is $C=(c_{ij})$, $i,j \in I = \{1, \dots, n\}$. 
Set $r^{\vee}$ to be the maximal number of edges connecting two vertices in the Dynkin diagram of $\lie g$ and define a renormalized form
$$(\cdot,\cdot) = r^{\vee}<\cdot,\cdot>$$
Then let 
\begin{equation}
d_i = \frac{(\alpha_i,\alpha_i)}{2}
\end{equation}
and set $D=\text{diag}(d_i)$, so that $B=DC$ is symmetric.

We consider $\hlie g$, the loop (or affine) algebra
associated to $\lie g$ (we do not consider central extension). 
The matrices $B,C,D$ can be extended to their affine counterparts
$\hat B,\hat C,\hat D$, of size $n+1$. Their entries will be
denoted by $b_{ij},c_{ij},d_i$ as before, but with 
$i,j$ running from $0$ to $n$. For $q \in \blb C^*$ (not a root of 1), the Hopf algebra $U_q(\hlie g)$ is the algebra generated by $k_i^{\pm 1}, x^{\pm}_i, i=0,\dots, n,$ with relations
\begin{gather} \notag
k_ik_i^{-1} = k_i^{-1}k_i = 1 \quad\quad k_ik_j = k_jk_i \quad\quad k_0\prod_{i=1}^{n}k_i^{\theta_i}=1\\ \notag
k_ix_j^{\pm}k_i^{-1} = q_i^{\pm c_{ij}}x_j^{\pm}\\
[x_i^+, x_j^+] = \delta_{ij} \frac{k_i - k_i^{-1}}{q_i - q_i^{-1}}\\ \notag
\sum_{m=0}^{1-c_{ij}} \begin{bmatrix} 1-c_{ij}\\ m \end{bmatrix}_{q_i} (x_i^{\pm})^m \,x_j^{\pm}\, (x_i^{\pm})^{1-c_{ij}-m} = 0, \quad i\neq j
\end{gather}
Here $q_i = q^{d_i}$ and {\scriptsize$\begin{bmatrix} r\\ m \end{bmatrix}_q$} $= \frac{[r]_q!}{[m]_q![r-m]_q!}$ where $[m]_q = \frac{q^m - q^{-m}}{q-q^{-1}}$ and $[m]_q! = [m]_q[m-1]_q\dots [1]_q$. The coalgebra structure and the antipode are given by
\begin{gather}\notag
\Delta(k_i^{\pm 1}) = k_i^{\pm 1}\otimes k_i^{\pm 1} \quad\quad \Delta(x_i^+) = x_i^+\otimes k_i + 1\otimes x_i^+ \quad\quad \Delta(x_i^-) = x_i^-\otimes 1 + k_i^{-1}\otimes x_i^-\\
\varepsilon(k_i) = 1 \quad\quad \varepsilon(x_i^{\pm}) = 0\\ \notag
S(k_i) = k_i^{-1} \quad\quad S(x_i^+) = -x_i^+ k_i^{-1} \quad\quad S(x_i^-) = -k_ix_i^- 
\end{gather}   
$U_q(\hlie g)$ has (in an appropriate sense \cite{EFK98}) a universal R-matrix that we denote by $\cal R$.


Given a representation $Y$ of $U_q(\hlie g)$ and $\lambda \in \lie h^*$, we set
$$Y[\lambda] = \{v\in Y; k_iv = q^{(\lambda,\alpha_i)}v\}$$ 
We will consider the category \cal C of finite dimensional (type {\bf 1}) representations of $U_q(\hlie g)$. It consists of finite dimensional representations with weight decomposition
$$Y = \bigoplus_{\lambda} Y[\lambda]$$
for $\lambda$ in the weight lattice of $\lie g$. 

Let $\text{\bf Gr}(\cal C)$ be the Grothendiek ring of $\mathcal C$. We will need the following theorem \cite{FreRe00}.

\begin{thm}\label{p:grring}
$\text{\bf Gr}(\cal C)$ is a commutative ring.
\end{thm}


\section{The Elliptic Central Character}

Recall that if $A$ is a finite dimensional algebra, then blocks 
in the category A-mod are parametrized by characters of the
center, $\chi: Z(A)\to \blb C$. 

Of course, the algebra
$U_q(\widehat{{\mathfrak g}})$ is not finite dimensional, and 
has a trivial center.  
Nevertheless, we will define a nontrivial analog of the notion 
of a central character (the {\it elliptic central character})
for $U_q(\widehat{{\mathfrak g}})$, 
which will allow us to compute the block decomposition
of the category of its finite dimensional representations.

For $z \in \blb C^*$, let $D_z$ be the automorphism of $U_q(\hlie g)$ given by $D_z(x_0^+) = zx_0^+, D_z(x_0^-)=z^{-1}x_0^-$ and the identity on the other generators. Given $X \in \cal C$ we can consider the family of shifted representations $X(z)$ obtained from $X$ by composing with $D_z$.

Let $X,Y \in \cal C$ and $\cal R|_{X(z)\otimes Y}$ be denoted by $R_{X,Y}(z)$. $R_{X,Y}(z)$ is a meromorphic function of $z$ regular at $0$ \cite{KaSo95,EtMo01}. Define
\begin{equation}
\eta_{X,Y}(z) = R^{21}_{Y,X}(z^{-1})R_{X,Y}(z)\in {\rm End}_{U_q(\hlie g)}(X(z)\otimes Y).
\end{equation}
Then $\eta_{X,Y}$ is an  elliptic function of $z$ with period
$q^{2r^{\vee}h^{\vee}}$ \cite{KaSo95}, i.e., a meromorphic function on the elliptic curve 
$$E = \frac{\blb C^*}{q^{2r^{\vee}h^{\vee}\blb Z}}$$
where $h^{\vee}$ is the dual Coxeter number of $\lie g$.

Let ${\rm Id}_E$ denote the identity functor 
of the category $\mathcal C_E:=\mathcal C\otimes_{\blb C}\blb
C(E)$.
 
\begin{prop}\label{endid}
If $X$ is irreducible, then there exists an element $\xi_X(z) \in
{\rm{End}}_{U_q(\hlie g)} ({\rm Id}_E)$, such that $\eta_{X,Y}(z) = 1\otimes\xi_X(z)|_Y$.
\end{prop}

Proposition \ref{endid} follows from the following lemma.

\begin{lem}\label{tensor} Let $Y\in \mathcal C$. Given a simple object $X\in \mathcal C$,
the map $\xi \mapsto 1\otimes\xi$, defines an
isomorphism  $\text{End}_{U_q(\hlie g)}(Y) \cong
\text{End}_{U_q(\hlie g)}(X(z)\otimes Y)$ for almost all $z\in
\blb C^*$.
\end{lem}

\begin{proof}
We have
$${\rm End}_{U_q(\hlie g)}(X(z)\otimes Y) \cong {\rm Hom}_{U_q(\hlie g)}(Y, \,^*X(z)\otimes X(z)\otimes Y) \cong {\rm Hom}_{U_q(\hlie g)}(Y\otimes\,^*Y,(\,^*X\otimes X)(z))$$
Let $Z_1, \dots, Z_n$ be the nontrivial constituents of a
composition series of $^*X\otimes X$. Then, for almost all $z$,
none of the $Z_i(z)$ occurs as a constituent in $Y\otimes\,^*Y$
(as $Z_i(z)$ are pairwise non-isomorphic for fixed $i$ and $z\in \blb C^*$)
and, consequently, the image of any morphism $f:Y\otimes\,^*Y \to
(^*X\otimes X)(z)$ has only trivial constituents. 

It is easy to show that ${\rm
  Ext}^1_{U_q(\hlie g)}(\blb C,\blb C) = 0$. Thus, the image of $f$ is
trivial, i.e., either zero or 1-dimensional (since $X$ is simple). The lemma is proved.
\end{proof}

\begin{cor}\hfill\\ \vspace{-.75cm}
\begin{enumerate}
\item If $Y$ is irreducible, then $\xi_X(z)|_Y$ is a scalar
  operator, 
and $\xi_X(z)|_Y=\xi_Y(z^{-1})|_X$.
\item $\xi_X(z)|_{Y_1\otimes Y_2} =
  \xi_X(z)|_{Y_1}\otimes\xi_X(z)|_{Y_2}$.
In particular, if $Y_i$ are irreducible, and $Y$ is a subquotient 
in $Y_1\otimes Y_2$, then
$\xi_X(z)|_Y=\xi_X(z)|_{Y_1}\xi_X(z)|_{Y_2}\in \Bbb C$. 

\item $\xi_X(z)|_{Y^*} = ((\xi_X(z)|_Y)^{-1})^*$
\item $\xi_X(z)|_{Y(u)} = \xi_X(\frac{z}{u})|_Y$ and $\xi_{X(w)}(z) = \xi_X(zw)$.
\item If an irreducible $X$ is a subquotioent of the tensor product of two irreducible $X_1,X_2$, then $\xi_X = \xi_{X_1}\xi_{X_2} = \xi_{X_2}\xi_{X_1}$.
\end{enumerate}
\end{cor}

\begin{proof}
The first statement immediately follows from Lemma \ref{tensor}. The
second and the last follow from the fusion laws for the
R-matrix. The third follows from the antipode property of the
$R$-matrix ($(S\otimes 1)(\mathcal R)=\mathcal R^{-1}$). The
fourth is a consequence of the fact that the same is true for $R_{X,Y}(z)$.
\end{proof}

\begin{defn}
Let $\cal I \subset {\rm Ob}(\cal C)$ be the set of isomorphism
classes of simple objects. An {\bf elliptic central character} in \cal C is an element $\chi \in \blb C(E)^{\cal I}$ such that 
\begin{enumerate}
\item $\chi_{X(w)}(z) = \chi_X(zw)$.
\item If $X$ is a subquotient of $X_1\otimes X_2$, for $X,X_i\in \cal I$, then $\chi_X = \chi_{X_1}\chi_{X_2}$.
\item $\chi_{\blb C} = 1$.
\end{enumerate}
\end{defn}

Given an elliptic central character $\chi$, set 
$$\cal C_{\chi} = \{Y\in \cal C; \xi_X(z)|_Z = \chi_X(z) \text{
  for all } Z\in \cal I, 
Z \text{ is a constituent of } Y\}$$
(it is possible that $\mathcal C_\chi=0$). If $Y\in \mathcal
C_\chi$, we will say that the elliptic central character of $Y$
is $\chi$. 

It immediately follows that \cal C is the direct sum of all $\cal
C_{\chi}$. 

We are ready to state our main result.

\begin{thm}\label{t:main}
The categories $\cal C_{\chi}$ are indecomposable.
In other words, they are the blocks of $\mathcal C$. 
\end{thm}

The set of $\chi$ for which $\cal C_\chi$ is nonzero can be
explicitly described. This description will be clear from the
proof. It will also be clear from the proof that 
if elliptic central characters of two representations coincide up
to scaling, then they coincide; so an elliptic central character 
of a representation is completely determined by its divisor of
zeros and poles on the elliptic curve $E$. 


\section{The $\lie{sl}_{n+1}$ Case}

To prove Theorem \ref{t:main} for all $\lie g$ we will need to use another realization of $U_q(\hlie g)$ using Drinfeld ``loop-like'' generators. However, due to the existence of Jimbo's algebra homomorphism $U_q(\hlie g) \to U_q(\lie g)$ when $\lie g$ is $\lie{sl}_{n+1}$ \cite{Jim86,EFK98}, we can already proceed with the proof in this case.

We have $r^\vee=1$, $h^\vee=n+1$ (but we will keep 
the notation $h^\vee$ to maintain similarities with types B-G). 
Let $V=\blb C^{n+1}$ be the vector representation of
$U_q(\hlie{sl}_{n+1})$. The following proposition is well known (see \cite{ChPr93,ChPr96} for example).

\begin{prop}\label{p:vectrep}
Any irreducible object of $\cal C$ is a subquotient of a tensor product of the form $V(z_1)\otimes \dots \otimes V(z_m)$. 
\end{prop}

\begin{cor}
Any elliptic central character $\chi$ is determined by $\chi_V$.
\end{cor}

The following is a consequence of the results in \cite{Kas00,Cha01}.

\begin{prop}\label{p:indec}
For a tensor product $V(z_1)\otimes \dots \otimes V(z_m)$ to be
cyclic on the highest weight vector (hence, indecomposable) 
it suffices that $z_j/z_k \neq q^{2}$ for $j<k$.
\end{prop}

\begin{prop}\label{p:iso} 
If $z_1, \dots , z_m$ is a sequence satisfying the condition in Proposition \ref{p:indec}, and $s \in S_m$ is such that $z_{s(1)}, \dots , z_{s(m)}$ also satisfies the condition, then the corresponding tensor products are isomorphic.
\end{prop}

\begin{proof} The $R$-matrix $\bar R_{V,V}(z)$ has singularities
  at $z=q^{\pm 2}$.   
Since for the transposed factors we have $z_j/z_k \neq q^{\pm
  2}$, the isomorphism 
is given by the action of $P\bar{R}_{V,V}(z_j/z_k)$, where
$\bar{R}$ is the 
normalized R-matrix of Proposition \ref{p:factor} below, and $P$ is the flip map.
\end{proof}

\begin{defn}\label{d:mono}
We say that a sequence $z_1, \dots, z_m$ is non-resonant if it
satisfies the condition of Proposition \ref{p:indec}. 
\end{defn}

Since any sequence can be arranged in a non-resonant order, we shall denote by $Y(z_1,\dots, z_m)$ any of the corresponding isomorphic indecomposable tensor products obtained from $z_1, \dots, z_m$.

It follows from Lemma \ref{l:linkcond} and Propositions \ref{p:vectrep}, \ref{p:indec} and \ref{p:iso}, that to prove Theorem \ref{t:main} for $U_q(\hlie{sl}_{n+1})$, it is enough to show that if $Y(z_1,\dots, z_m)$ and $Y(w_1, \dots, w_k)$ have the same elliptic central character, then they are linked.

\begin{lem}\label{l:linked}
For any $w\in \blb C^*$, $Y(z_1,\dots, z_m)$ is linked to $Y(z_1,\dots, z_m, w, q^{2}w,\dots, q^{2(h^{\vee}-1)}w)$. In particular $Y(z_1,\dots, z_m)$ is linked to $Y(z_1,\dots, z_{j-1}, z_jq^{2h^{\vee}}, z_{j+1},\dots, z_m)$.
\end{lem}

\begin{proof}
Since the trivial representation is contained in
$Y(w,q^2w,\dots,q^{2(h^{\vee}-1)}w)$ 
for any $w$ (as ``the top quantum exterior power'' of $V(w)$), it follows (using the fact that $\text{\bf Gr}(\cal C)$ is commutative) that any simple constituent of $Y(z_1,\dots, z_m)$ is also a constituent of $Y(z_1,\dots, z_m, w, q^{2}w,\dots, q^{2(h^{\vee}-1)}w)$. For the second statement, let $Y_1, Y_2$ be simple constituents in each of the considered tensor products respectively.  Then both $Y_i$ are subquotients of $Y(z_1,\dots,z_j,q^2z_j,\dots,q^{2h^{\vee}}z_j, z_{j+1},\dots, z_m)$.
\end{proof}

In light of Lemma \ref{l:linked} it remains to show that if $Y(z_1,\dots, z_m)$
and $Y(w_1,\dots, w_k)$ have the same elliptic central character,
then the sequence $(z_1,...,z_m)$ can be obtained from 
$(w_1,...,w_k)$ by 
permutations and by adding and removing sequences $(z,q^2z,...,q^{2(h^\vee-1)}z)$
(which includes the transformations $z_j\to z_jq^{2h^\vee}$).   

To do this, we will write down the formula for $\xi_V(z)|_V$ and
analyze its singularity structure. Drinfeld realization will be
used to do this in the other cases, but for $\lie{sl}_{n+1}$ one
can compute it in a more ``naive'' way, as was (essentially) done
in \cite{Mou01}.
We have

\begin{equation}
\xi_V(z)|_V = q^{\frac{2(h^{\vee}-1)}{h^{\vee}}}\prod_{j=0}^{\infty} \varrho(q^{2jh^{\vee}}z)\varrho(q^{2jh^{\vee}}z^{-1})
\end{equation}
where
\begin{equation}\label{eq:fro}
\varrho(z) = \frac{(1-z)(1-zq^{2h^{\vee}})}{(1-zq^2)(1-zq^{2(h^{\vee}-1)})}
\end{equation} 
Then, the structure of zeros and poles of $\xi_V(z)|_V$ on $E$ is given by the following pictures
$$\begin{matrix}
\bullet^1 & \bullet^{q^2} & \bullet^{q^4} & \dots & \bullet^{q^{2(h^{\vee}-2)}}    &\bullet^{q^{2(h^{\vee}-1)}}\\
\hspace{-.15cm}2 & \hspace{-.35cm}-1 & \hspace{-.25cm}0 &\dots &\hspace{-1.1cm}0          &\hspace{-1.1cm}-1
\end{matrix} \qquad\text{for }n\geq 2\qquad\text{and}\quad\quad
\begin{matrix}
\bullet^1 & \bullet^{q^2}\\
\hspace{-.15cm}2 & \hspace{-.35cm}-2
\end{matrix} \qquad\text{for }n=1$$
where positive numbers stand for zeros (of that order) and
negative for poles. 

The fact that the trivial representation is contained in
$Y(w,q^2w,\dots,q^{2(h^{\vee}-1)}w)$ is reflected in the relation 
\begin{equation}\label{rela}
\prod_{s=0}^{h^\vee-1}\xi_V(zw^{-1}q^{-2s})=1. 
\end{equation}
To prove our claim, we need to show that any multiplicative
relation between $\xi_V(zu)$, $u\in \blb C^*$ is a combination
of relations of the form (\ref{rela}). For this, it suffices to
show that the functions
$\xi_V(z), \xi_V(zq^{-2})$,...,$\xi_V(zq^{-2(h^\vee-2)})$ are
multiplicatively independent (for $\lie{sl}_2$ this is clear). 

To do this, we will rephrase the problem in a linear algebra
setting. Consider the group $\blb Z^{h^{\vee}-1}$. To the
function 
$\xi_V(zq^{-2s})|_{V}$, $0\le s\le h^\vee-2$ 
assign a vector $v_s$ in $\blb Z^{h^{\vee}-1}$
given by
\begin{alignat*}3
& v_0 &=& \,(2,-1, 0, 0, \dots, 0)\\
& v_{h^{\vee}-2} &=& \,(0, 0, \dots, 0, -1, 2)\\
& v_s &=& \,(0,\dots, 0, -1,2,-1, 0,\dots, 0) \qquad\text{for } 0< s <h^{\vee}-2
\end{alignat*}
where the $2$ is the $s$-th entry, if we label them from $0$ to
$h^{\vee}-2$. 
The entries of these vectors are the orders of
 the singularities of
$\xi_V(zq^{-2s})$ on the sequence $1, q^2,..., q^{2(h^\vee-2)}$. 
Then we are left to show that the vectors $v_s$ are linearly
independent. This is equivalent to showing that the matrix
$T_n$, whose rows are the vectors $v_s$, has a nonvanishing
determinant. But one easily sees that $\det T_n = 2\det T_{n-1}
-\det T_{n-2}$ and use induction to get $\det T_n = n+1$
(in fact, $T_n$ is the Cartan matrix of type $A_n$). 

This proves Theorem \ref{t:main} for $\lie{sl}_{n+1}$.


\section{Drinfeld Realization}

To compute the elliptic central characters for the other $\lie g$
we will need to use a suitable formula for $\cal R$. This formula
was found in \cite{KhTo92} and we will use the version
\cite{FreRe00} 
involving the Drinfeld realization \cite{Dri88} of $U_q(\hlie g)$ in terms of ``loop like'' generators. 
In fact, we will use the following (slightly different) realization proposed in \cite{Bec93}.

\begin{thm}
$U_q(\hlie g)$ is isomorphic to the algebra with generators $k_i^{\pm 1}, \,h_{i,l}, \, x_{i,r}^{\pm}$, 
where $i\in I, \,l\in \blb Z\backslash \{0\}$ and $r \in \blb Z$ with defining relations
\begin{gather}\notag
k_ik_i^{-1} = k_i^{-1}k_i = 1\\ \notag
k_ik_j = k_jk_i \quad\quad k_i h_{j,r} = h_{j,r}k_i\\ \notag
k_ix_{j,r}{\pm}k_i^{-1} = q_i^{\pm b_{ij}}x_{j,r}^{\pm} \quad\quad 
[h_{i,r}, x_{j,s}^{\pm}] = \pm\frac{1}{r}[rc_{ij}]_{q^i}\, x_{j,r+s}\\ 
x_{i,r+1}^{\pm}x_{j,s}^{\pm} - q_i^{\pm c_{ij}}x_{j,s}^{\pm}x_{i,r+1}^{\pm} = q_i^{\pm c_{ij}} x_{i,r}^{\pm}x_{j,s+1}^{\pm} - x_{j,s+1}^{\pm}x_{i,r}\\ \notag
[x_{i,r}^+,x_{j,s}^-]= \frac{\delta_{ij}}{q_i - q_i^{-1}}\, \big( \phi_{i,r+s}^+ - \phi_{i,r+s}^-\big)\\ \notag
\sum_{s\in S_m}\sum_{k=0}^m (-1)^k\begin{bmatrix} m\\ k\end{bmatrix}_{q_i} x_{i,r_{s(1)}}^{\pm}\dots x_{i,r_{s(k)}}^{\pm} x_{j,s}^{\pm}\, x_{i,r_{s(k+1)}}^{\pm}\dots x_{i,r_{s(m)}}^{\pm} =0\quad \text{if }i\neq j
\end{gather}
where $r_1, \dots, r_m \in \blb Z$, $m=1-c_{ij}$, $S_m$ is the symmetric group on $m$ symbols and $\phi_{i,r}^{\pm}$ are given by the identity of power series
\begin{equation}
\sum_{r=0}^{\infty} \phi_{i,r}^{\pm} u^{\pm r} = k_i^{\pm 1}\exp\Big( \pm(q - q^{-1})\sum_{s=1}^{\infty} h_{i,\pm s}u^{\pm s}\Big)
\end{equation}
\end{thm} 

\obs In fact this version follows the notation in
\cite{FreRe00}. But observe that the comultiplication in
\cite{FreRe00} differs from ours. They are connected by the automorphism of $U_q(\hlie g)$ sending $k_i$ to $k_i^{-1}$, $q$ to $q^{-1}$ and keeping fixed the other generators. This will give rise to some differences in the signs.
We use the original definitions of \cite{Dri85,Jim85} which also coincides with \cite{ChPr94bo,ChPr95} and \cite{Mou01, EtMo01} (these last two are the reasons for our choice). But we remark that the definition in \cite{FreRe00} has been used also in \cite{Cha01,FreMu01} and many other most recent works.

Using the ``loop like'' generators, to each fundamental weight
$\omega_i$, 
one may associate a family of (shifted) {\it fundamental representations} $V_i(z)$. In general, the fundamental representations are not irreducible as $U_q(\lie g)$-modules, but we still have the following theorem \cite{ChPr94bo,ChPr95}.

\begin{thm}\label{t:fundrep}
Any irreducible representation of $U_q(\hlie g)$ is isomorphic to a subquotient of a tensor product of shifted fundamental representations.
\end{thm}

\obs Using the Drinfeld realization, Chari and Pressley (see e.g.\cite{ChPr95}) defined the concept of an affinization of 
$V_{\lambda}$ for any finite dimensional $U_q(\lie g)$-module $V_{\lambda}$ with highest weight $\lambda$. 
By definition, an irreducible affinization of $V_{\lambda}$ is a $U_q(\hlie g)$-representation  $V$ isomorphic to one of the form 
$V_{\lambda} \oplus\bigoplus_{\mu < \lambda} V_{\mu}$ as a representation of $U_q(\lie g)$. 
We refer to \cite{ChPr95} for more details.
We also mention one fact that will be used later. 
Namely, if $\lambda = \sum \lambda_i\omega_i$, then all irreducible affinizations of $V_{\lambda}$ are obtained as subquotients of  tensor products
of the form $\otimes_i^n(\otimes_{j=1}^{\lambda_i}V_i(z_{ij}))$.

The following proposition \cite{FreR, EFK98} 
is our first tool to compute $\xi_X(z)|_Y$.

\begin{prop}\label{p:factor}
Let $X,Y$ be irreducible representations of $U_q(\hlie g)$. Then
\begin{equation}
R_{X,Y}(z) = f_{X,Y}(z)\bar{R}_{X,Y}(z)
\end{equation}
where $f_{X,Y}$ is a scalar meromorphic function in \blb C, regular at $0$ with $f_{X,Y}(0) \neq 0$, and the matrix elements of $\bar{R}_{X,Y}(z)$ are rational functions of $z$ regular at $0$ and such that $\bar{R}_{X,Y}(z)(x_0\otimes y_0) = x_0\otimes y_0$. Here $x_0, y_0$ are the highest weight vectors of $X$ and $Y$ as $U_q(\lie g)$-modules. If $|q|<1$, $f_{X,Y}$ can be represented as
$$f_{X,Y}(z) = q^{(\lambda,\mu)}\prod_{j=0}^{\infty} \varrho_{X,Y}(q^{2jr^{\vee}h^{\vee}}z)$$
where $\lambda$ and $\mu$ are the highest weights of $X$ and $Y$ and $\varrho_{X,Y}$ is a rational function such that $\varrho_{X,Y}(0)=1$. Furthermore, $\bar{R}$ is unitary:
\begin{equation}\label{eq:unit}
\bar{R}_{Y,X}^{21}(z^{-1})\bar{R}_{X,Y}(z) = 1
\end{equation}
\end{prop}

\begin{cor}
If $X,Y$ are irreducible, $\xi_X(z)|_Y$ is the scalar operator given by $f_{X,Y}(z)f_{X,Y}(z^{-1})$.
\end{cor}

It was proved in \cite{FreMu01} that if $\bar{R}_{V_i,V_j}(z)$ is
not regular at $z_0$, then $z_0$ must belong to the set $\cal P = \{q^k; 2\leq k \leq
r^{\vee}h^{\vee}, k\in\blb Z\}$ and, if $\bar{R}_{V_i,V_j}(z_0)$
is not invertible, then $z_0 \in \cal P^{-1} = \{q^{-k}; 2\leq k
\leq r^{\vee}h^{\vee},k\in \blb Z\}$.

\begin{cor}\label{c:poles}
Let $\varrho_{V_i, V_j}(z)$ be the function of Proposition \ref{p:factor} and $\cal P_{ij}$ be the subset of \cal P where $\varrho_{V_i, V_j}(z)$ has a pole. Then $\bar{R}_{V_i,V_j}(z)$ is not invertible exactly on $\cal P_{ij}^{-1}$.
\end{cor}

\begin{proof}
The function $\varrho_{X,Y}(z)$ is characterized by the following equation \cite{EFK98}:
\begin{equation}
((\bar{R}_{X,Y}(z))^{-1})^{t_1} = \varrho_{X,Y}(z) ((\bar{R}_{X,Y}(q^{2r^{\vee}h^{\vee}}z))^{t_1})^{-1}
\end{equation}
where $(\sum a_j \otimes b_j)^{t_1} = \sum a_j^* \otimes
b_j$. Recall that $V_j(z)^* \cong V_{j^*}(q^{r^{\vee}h^{\vee}})$
where $V_{j^*}$ is the fundamental representation of  $U_q(\lie
g)$ dual to $V_j$. Then, if $z_0 \in P_{ij}^{\pm 1}$ and $\bar{R}_{V_i,V_j}(z_0)$ is
invertible, both $((\bar{R}_{V_i,V_j}(z))^{-1})^{t_1}$ and
$((\bar{R}_{V_i,V_j}(q^{2r^{\vee}h^{\vee}}z))^{t_1})^{\pm 1}$ are
regular at $z_0$. Hence so is $\varrho_{V_i,V_j}(z)$. Conversely,
if $\bar{R}_{V_i,V_j}(z_0)$ is not invertible, then
$((\bar{R}_{V_i,V_j}(z))^{-1})^{t_1}$ has a pole at $z_0$, but
$((\bar{R}_{V_i,V_j}(q^{2r^{\vee}h^{\vee}}z))^{t_1})^{-1}$ is
still regular. Thus, $\rho_{X,Y}$ must have a pole at $z_0$.  
\end{proof}

The following corollary will be useful later. 
Let $\cal S_{ij} = \cal P_{ij}\cup\cal P_{ij}^{-1}$.

\begin{cor}\label{c:smallshift}
Suppose that $V_r(q^p)$ is a subrepresentation of $V_i\otimes V_j(q^l)$. 
Then, for any $m =1,\dots, n$, we have
$\cal S_{mr}q^p\subset \cal S_{mi}\cup\cal S_{mj}q^l$. 
\end{cor}

\begin{proof}
Given $m$ and $z\in \blb C$, consider the inclusion
$$V_m(z)\otimes V_r(q^p) \hookrightarrow V_m(z)\otimes V_i\otimes V_j(q^l)$$
By the fusion laws for the universal R-matrix, the singularities of 
$\bar R_{V_m,V_i\otimes V_j(q^l)}(z)$ (poles and points where it is not invertible),
must be contained in $\cal S_{mi}\cup\cal S_{mj}q^l$. On the other hand, let $\bar R_{m,ij}^r(z)$ denote the restriction of 
$\bar R_{V_m,V_i\otimes V_j(q^l)}(z)$ to  
$V_m(z)\otimes V_r(q^p)$. Then $g(z) \bar R_{m,ij}^r(z)=
\bar R_{V_m,V_r(q^p)}(z)$, for some rational function $g(z)$. 
Suppose that, at $z_0$, $\bar R_{V_m,V_r(q^p)}(z)$ is not invertible,
but $\bar R_{m,ij}^r(z_0)$ is defined and invertible. We conclude that $g(z_0)=0$. But this would imply that 
$\bar R_{V_m,V_r(q^p)}(z_0)=0$, what is impossible by the normalization of $\bar R_{V_m,V_r(q^p)}(z)$. Since is $z \mapsto z^{-1}$ 
is a bijection $\cal P_{mr} \to \cal P_{mr}^{-1}$, we conclude that $\bar R_{m,ij}^r(z_0)$ has a singularity whenever
$\bar R_{V_m,V_r(q^p)}(z)$ has a pole.
\end{proof}  

Now we need a tool to calculate $f_{X,Y}$ for fundamental $X$ and $Y$. It is the formula for \cal R found in \cite{KhTo92}. The method we will describe now was developed for general irreducibles $X,Y$ in \cite{FreRe00}. So our computation is a specialization of those in \cite{FreRe00}. First we define the matrices $B(q), D(q)$ and $M(q)$ to be $b_{ij}(q) = [b_{ij}]_q, \,d_{ij}(q) = \delta_{ij}[d_i]_q$, for $i,j = {1, \dots, n}$ and $M(q) = D(q)\tilde{B}(q)D(q)$ where $\tilde{B}(q) = B(q)^{-1}$.

\begin{thm}
\cite{KhTo92}
The universal R-matrix \cal R of $U_q(\hlie g)$ can be
represented in the form
\begin{equation}
\cal R = \cal R^+\breve{\mathcal{R}}\cal R^-\cal R^0
\end{equation}
where $\cal R^{\pm} \in U_q(\hlie n_{\pm})\otimes U_q(\hlie n_{\mp})$, $\cal R^0(x\otimes y) = q^{(\lambda,\mu)}x\otimes y$ if $x,y$ have weight $\lambda$ and $\mu$ respectively, and the ``imaginary'' R-matrix $\breve{\mathcal{R}}$ is given by
\begin{equation}
\breve{\mathcal{R}} = \exp\Big((q-q^{-1})\sum_{k>0, i,j} \frac{k}{[k]_q} \tilde{b}_{ij}(q^k)h_{i,k}\otimes h_{j,-k}\Big)
\end{equation}
\end{thm}

Then we can compute $f_{X,Y}$ by calculating the action of \cal R on the tensor product of the corresponding highest weight vectors. $\cal R^{\pm}$ will act as the identity while $\cal R^0$ will contribute a constant. Hence, the essential information is contained in $\breve{\cal R}$. Let $v_i$ be the highest weight vector of $V_i(z)$. The action of $h_{j,k}$ on $v_i$ is given by \cite{FreRe00}
\begin{equation}
h_{j,k}v_i = \delta_{ij} \frac{(q_j^k-q_j^{-k})z^k}{(q-q^{-1})k}\,v_j
\end{equation}
Consequently
\begin{equation}\label{e:sing}
\breve{\cal R}(v_i\otimes v_j) = \exp\Big(\sum_{k>0}m_{ij}(q^k)(q^k-q^{-k})\frac{z^k}{k}\Big)v_i\otimes v_j = q^{-(\omega_i,\omega_j)}f_{V_i,V_j}(z)v_i\otimes v_j
\end{equation}
The singularity structure for $f_{V_i,V_j}(z)$ can then be read from the matrix
$M(q)$ in the following way.
The term $m_{ij}(q^k)(q^k-q^{-k})\frac{z^k}{k}$ will be of the form
$\frac{\pi(q^k)}{1-q^{pk}}\frac{z^k}{k}$,
where $\pi(q)$ is a Laurent polynomial in $q$
(symmetric under $q \mapsto q^{-1}$) and $p=2r^{\vee}h^{\vee}$.
Then, each monomial $\pm q^m$ in $\pi(q)$
will contribute with a factor
$\prod_{l\geq 0}(1-q^mzq^{pl})^{\mp1}$ to $f_{V_i,V_j}(z)$
(recalling that $\exp(-\sum_k \frac{y^k}{k}) = \exp(\log(1-y)) = 1-y$
and that $\frac{1}{1-x}=\sum_l x^l$).

\obs The matrices $M(q)$ are listed in the appendix. For the
classical algebras they were listed in \cite{FreRe98} (with a few misprints). For types $E$ and $F$ we used the software Mathematica to compute $M(q)$.


\section{Proof of Theorem \ref{t:main} : The Remaining Cases}

In this section we will use the notation $\xi_{ij}(z) = \xi_{V_i}(z)|_{V_j}$.

We state the following version of Proposition \ref{p:vectrep}.

\begin{thm}\label{t:generator}
Every irreducible object of \cal C is a subquotient of a tensor product of the form $V_{i_1}(z_1)\otimes \dots\otimes V_{i_m}(z_m)$ where $i_j$ run through the indices of the black nodes of the Dynkin diagram of $\lie g$ in table 1.
\end{thm}

For classical $\lie g$ it is proved in \cite{ChPr96}. For $\lie
e_6$ it is deduced immediately from the computations in
\cite{ChPr91}. We will prove it in the remaining cases together
with the proof of Theorem \ref{t:main}. (We note that as far as
we know, in the
cases $F_4$ and $G_2$, Theorem \ref{t:generator} was also proved
by Chari and Pressley in 1991 when preparing the paper 
\cite{ChPr91}; however, the proof is not written in the paper,
since the result was not of interest at that time).  

\pagebreak
\vspace{.5cm}

{\centerline {\bf Table 1}}
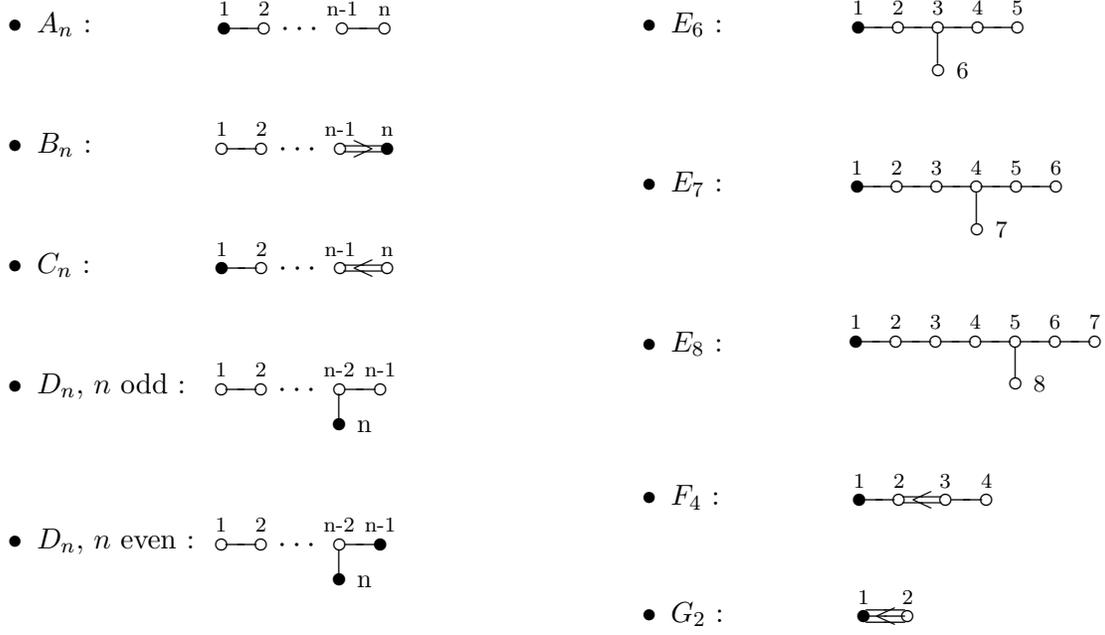
\vspace{.5cm}
\begin{multicols}{2}
\begin{itemize}

\item $A_n$ : {\large \vspace{-.55cm} $$\stackrel{1}{\bullet}\hspace{-.2cm}\sn\hspace{-.2cm}\stackrel{2}{\circ} \dots \stackrel{\text{n-1}}{\circ}\hspace{-.29cm}\sn\hspace{-.18cm}\stackrel{\text{n}}{\circ}$$}

\item $B_n$ : {\large \vspace{-.55cm} $$\stackrel{1}{\circ}\hspace{-.2cm}\sn\hspace{-.2cm}\stackrel{2}{\circ} \dots \stackrel{\text{n-1}}{\circ}\hspace{-.18cm}\dnl\hspace{-.07cm}\stackrel{\text{n}}{\bullet}$$}

\item $C_n$ : {\large \vspace{-.55cm} $$\stackrel{1}{\bullet}\hspace{-.2cm}\sn\hspace{-.2cm}\stackrel{2}{\circ} \dots \stackrel{\text{n-1}}{\circ}\hspace{-.18cm}\dnr\hspace{-.07cm}\stackrel{\text{n}}{\circ}$$}

\item $D_n$, $n$ odd : {\large \vspace{-.95cm} $$\stackrel{1}{\circ}\hspace{-.2cm}\sn\hspace{-.2cm}\stackrel{2}{\circ} \dots \stackrel{\text{n-2}}{\circ}\hspace{-.29cm}\sn\hspace{-.3cm}\stackrel{\text{n-1}}{\circ}$$
\vspace{-1.12cm}$$\hspace{.92cm}|$$
\vspace{-1.05cm}$$\hspace{1.25cm}\bullet\text{ \small n}$$}

\item $D_n$, $n$ even : {\large \vspace{-.95cm} $$\stackrel{1}{\circ}\hspace{-.2cm}\sn\hspace{-.2cm}\stackrel{2}{\circ} \dots \stackrel{\text{n-2}}{\circ}\hspace{-.29cm}\sn\hspace{-.3cm}\stackrel{\text{n-1}}{\bullet}$$
\vspace{-1.12cm}$$\hspace{.92cm}|$$
\vspace{-1.05cm}$$\hspace{1.25cm}\bullet\text{ \small n}$$}
\\
\item $E_6$  : {\large \vspace{-.55cm} $$\stackrel{1}{\bullet}\hspace{-.2cm}\sn\hspace{-.2cm}\stackrel{2}{\circ} \hspace{-.2cm}\sn\hspace{-.2cm} \stackrel{\text{3}}{\circ}\hspace{-.2cm}\sn\hspace{-.2cm}\stackrel{4}{\circ}\hspace{-.2cm}\sn \hspace{-.2cm}\stackrel{5}{\circ}$$
\vspace{-.93cm}$$\hspace{-.02cm}|$$
\vspace{-.93cm}$$\hspace{.3cm}\circ\text{ \footnotesize$6$}$$}

\item $E_7$  : {\large \vspace{-.55cm} $$\hspace{.5cm}\stackrel{1}{\bullet}\hspace{-.2cm}\sn\hspace{-.2cm}\stackrel{2}{\circ} \hspace{-.2cm}\sn\hspace{-.2cm} \stackrel{\text{3}}{\circ}\hspace{-.2cm}\sn\hspace{-.2cm}\stackrel{4}{\circ}\hspace{-.2cm}\sn \hspace{-.2cm}\stackrel{5}{\circ}\hspace{-.2cm}\sn\hspace{-.2cm}\stackrel{6}{\circ}$$
\vspace{-.93cm}$$\hspace{1.02cm}|$$
\vspace{-.93cm}$$\hspace{1.33cm}\circ\text{ \footnotesize$7$}$$}

\item $E_8$  : {\large \vspace{-1.cm} $$\hspace{1cm}\stackrel{1}{\bullet}\hspace{-.2cm}\sn\hspace{-.2cm}\stackrel{2}{\circ} \hspace{-.2cm}\sn\hspace{-.2cm} \stackrel{\text{3}}{\circ}\hspace{-.2cm}\sn\hspace{-.2cm}\stackrel{4}{\circ}\hspace{-.2cm}\sn \hspace{-.2cm}\stackrel{5}{\circ}\hspace{-.2cm}\sn\hspace{-.2cm}\stackrel{6}{\circ} \hspace{-.2cm}\sn\hspace{-.2cm}\stackrel{7}{\circ}$$
\vspace{-1.09cm}$$\hspace{2.05cm}|$$
\vspace{-.94cm}$$\hspace{2.35cm}\circ\text{ \footnotesize$8$}$$} 

\item $F_4$ : {\large \vspace{-.55cm} $$\hspace{-.4cm}\stackrel{1}{\bullet}\hspace{-.2cm}\sn\hspace{-.2cm}\stackrel{2}{\circ} \hspace{-.08cm}\dnr\hspace{-.08cm}\stackrel{3}{\circ} \hspace{-.18cm}\sn\hspace{-.2cm}\stackrel{4}{\circ}$$}

\item $G_2$ : {\large \vspace{-.55cm} $$\hspace{-1.4cm}\stackrel{1}{\bullet}\hspace{-.1cm}\equiv\hspace{-.2cm}<\hspace{-.2cm}\equiv \hspace{-.1cm}\stackrel{2}{\circ}$$}

\end{itemize}

\end{multicols}            

\begin{cor}
If $\lie g$ is not of type $D_n$ for $n$ even, then any elliptic
central character $\chi$ is determined by $\chi_{V}$, where $V=V_b$ and $b$ is the index of the black node in table 1. For $D_{n}$, when $n$ is even, $\chi$ is determined by its value on the two half spin representations $V_{n-1}$ and $V_n$.
\end{cor}

Combining the results in \cite{Kas00,Cha01} with Corollary \ref{c:poles}, we get the following stronger version of Proposition \ref{p:indec}, which is crucial in our proof of Theorem \ref{t:main}.

\begin{thm}
For a tensor product of fundamental representations
$V_{k_1}(z_1)\otimes \dots\otimes V_{k_l}(z_l)$ to be cyclic on
the highest weight vector (hence indecomposable), it suffices that
$$\frac{z_r}{z_s}\neq q^{2d_{k_s}+p} \quad\text{for}\quad r<s\, , \quad p \geq 0 \quad\text{and}\quad 2\leq 2k_s + p\leq r^{\vee}h^{\vee}$$
In other words, it suffices that $\bar R_{V_{k_r}, V_{k_s}}(\frac{z_r}{z_s})$ be regular for $r<s$.
\end{thm}

Then we can prove the corresponding version of Proposition
\ref{p:iso} in a similar way and define $Y(z_1,\dots, z_m)$
analogously to Definition \ref{d:mono}. For $D_n$, when $n$ is
even, we let the half spin representations be denoted by
$V_{\pm}$ and define $Y_{+}(z_1,\dots, z_m)$ to be any of the
isomorphic indecomposable tensor products obtained from
$z_1,\dots,z_m$ using only shifts of $V_{+}$.
Similarly, we define $Y_-(w_1,\dots, w_l)$ using $V_-$. 
Then we can define $Y(z_1,\dots, z_m|
w_1, \dots, w_l)$ to be any non-resonant 
(i.e., satisfying the cyclicity condition on the highest weight
vector) reordering of 
$Y_{+}(z_1,\dots, z_m)\otimes Y_-(w_1,\dots, w_l)$. 

Therefore, once we have computed the singularity structure of $\xi_{ij}$, the proof goes,
case by case, in a similar way it did for type $A_n$. Namely, the
proof consists of two steps. 

\noindent{\bf Step 1}. We prove a version of lemma \ref{l:linked}
stating the basic linking relations.

\noindent{\bf Step 2}. To conclude that, if $Y(z_1,\dots, z_m)$ and $Y(w_1,\dots, w_l)$ have
the same elliptic central character, then they are linked, we find all multiplicative 
relations between the considered $\xi_{ij}(zq^{-2s})$ 
(by computing the kernel of a certain integer matrix over $\Bbb
Z$) and check that they correspond to 
the linking relations of the lemma.


\subsection{Type $B_n$}

For $B_n$ we have $r^{\vee} = 2, h^{\vee} = 2n-1$. The black node corresponds to the spin representation $V_n$. The singularity arrangement for $\xi_{nn}(z)$ is
\begin{gather*}
\begin{matrix}
\bullet^{1} & \bullet^{q^2} & \bullet^{q^4} & \dots & \bullet^{q^{2h^{\vee}-2}}&  \bullet^{q^{2h^{\vee}}} &  \bullet^{q^{2h^{\vee}+2}}& \dots & \bullet^{q^{4h^{\vee}-4}} & \bullet^{q^{4h^{\vee}-2}}\\ 
\hspace{-.25cm}2 & \hspace{-.25cm}-1 & \hspace{-.25cm}1 & \dots & \hspace{-1cm}1 & \hspace{-.75cm}-2 & \hspace{-1cm}1 & \dots & \hspace{-1cm}1 & \hspace{-1.25cm}-1
\end{matrix}
\end{gather*}
where the dots mean that the sequence goes on like
\begin{gather*}
\begin{matrix}
\dots & \bullet^{q^{2(k-2)}} & \bullet^{q^{2(k-1)}} & \bullet^{q^{2k}} & \bullet^{q^{2(k+1)}} & \dots \\ 
\dots & \hspace{-.85cm}1 & \hspace{-1cm}-1 & \hspace{-.5cm}1 &  \hspace{-1cm}-1 &\dots
\end{matrix}
\end{gather*}
except at $1$ and $q^{2h^{\vee}}$.

The corresponding version of Lemma \ref{l:linked} is

\begin{lem}\label{l:linked2}
For any $w\in \blb C^*$, $Y(z_1,\dots, z_m)$ is linked to $Y(z_1,\dots, z_m, w,  q^{r^{\vee}h^{\vee}}w)$. In particular  $Y(z_1,\dots, z_m)$ is linked to $Y(z_1,\dots, z_{j-1}, z_jq^{2r^{\vee}h^{\vee}}, z_{j+1},\dots, z_m)$.
\end{lem}

\begin{proof}
It is proved analogously to Lemma \ref{l:linked} and follows from the fact that $\blb C$ occurs as a constituent of $Y(w,q^{r^{\vee}h^{\vee}}w) = V(w)\otimes V(q^{r^{\vee}h^{\vee}}w)$.
\end{proof} 

For the last part of the proof (i.e. the proof that there is no
relations between $\xi_{nn}(zu)$ other than given by Lemma
\ref{l:linked2}), 
we consider the assignment
$$\xi_{nn}(zq^{-2s}) \mapsto v_s
:=((-1)^s,(-1)^{s+1},\dots,-1,2,-1,\dots,(-1)^s)\in \blb Z^{h^\vee}$$
for $s=0, \dots, h^{\vee}-1$.
The entries of these vectors are the orders of the singularities of
$\xi_{nn}(zq^{-2s})$ on the sequence $1, q^2,..., q^{2h^\vee-2}$.

\obs Observe that these vectors contain information only about the ``first half'' of the 
singularity structure. Since $V_n$ is ``self dual'', the ``second half''
is obtained from the first one by a change of signs.

The corresponding matrix $T_n$ is the $h^{\vee}\times h^{\vee}$-matrix of the form
$$T_n = 
\begin{pmatrix}
b      & a      & -a     & \dots  & a      & -a\\
a      & b      & a      & \ddots & -a     & a\\
-a     & \ddots & \ddots & \ddots & \ddots & \vdots\\
\vdots & \ddots & \ddots & \ddots & \ddots & \vdots\\
a      & -a     & \ddots & a      & b      & a\\
-a     & a      & \dots  & -a     & a      & b
\end{pmatrix}$$
where $a=-1$ and $b=2$. Then $\det T_n =
(a+b)^{h^{\vee}-1}(b-(h^{\vee}-1)a) = 
h^{\vee}+1 = 2n$. Thus, the rows of this matrix are linearly
independent, and hence there is no additional relations, as
desired. 


\subsection{Type $C_n$}

Here $r^{\vee} = 2, h^{\vee}=n+1$ and, for $1\leq i\leq j\leq n$, the zeros and poles arrangement of $\xi_{ij}(z)$ is then given by
$$\begin{matrix}
\bullet^{q^{|j-i|}} & \bullet^{q^{j+i}} & \bullet^{q^{2h^{\vee}-(j+i)}} & \bullet^{q^{2h^{\vee}-|j-i|}} & \bullet^{q^{2h^{\vee}+|j-i|}} & \bullet^{q^{2h^{\vee}+j+i}} & \bullet^{q^{4h^{\vee}-(j+i)}} & \bullet^{q^{4h^{\vee}-|j-i|}}\\
\hspace{-.75cm}1 & \hspace{-.75cm}-1 & \hspace{-1.5cm}1 & \hspace{-1.5cm}-1 & \hspace{-1.5cm}-1 & \hspace{-1.25cm}1 & \hspace{-1.45cm}-1 & \hspace{-1.4cm}1 
\end{matrix}$$

The black node in table 1 corresponds to the natural vector representation $V=V_1$. Therefore we can restrict ourselves to analysing the singularities of $\xi_{11}$  
$$\begin{matrix}
\bullet^{1} & \bullet^{q^2} & \hspace{1cm}\bullet^{q^{2h^{\vee}-2}} & \bullet^{q^{2h^{\vee}}} &  \bullet^{q^{2h^{\vee}+2}} & \hspace{1cm}\bullet^{q^{4h^{\vee}-2}}\\
\hspace{-.25cm}2 & \hspace{-.25cm}-1 & \hspace{1cm}\hspace{-.75cm}1 & \hspace{-1.0cm}-2 & \hspace{-1.0cm}1 & \hspace{1cm}\hspace{-1.0cm}-1
\end{matrix}$$

Lemma \ref{l:linked2} remains valid as in 
the $B_n$ case. To perform the second step of the proof, we assign
\begin{alignat*}3
& \xi_{11}(z) & \quad\mapsto & \quad v_0 :=(2,-1, 0, 0, \dots, 0,0,1)\\
& \xi_{11}(zq^{-2}) & \quad\mapsto & \quad v_1 :=(-1,2,-1, 0, 0,\dots,0, 0)\\
& \dots\dots\dots & & \quad\dots\dots\dots\dots\dots\dots\dots\dots\\
& \xi_{11}(zq^{-2(h^{\vee}-2)}) & \quad\mapsto & \quad v_{h^{\vee}-2} :=(0, 0,\dots,0, 0, -1, 2, -1)\\
& \xi_{11}(zq^{-2(h^{\vee}-1)}) & \quad\mapsto & \quad v_{h^{\vee}-1} :=(1, 0, 0, \dots, 0, 0, -1, 2)\\
\end{alignat*}
in $\blb Z^{h^{\vee}}$. Again, the entries of these vectors are the orders of the singularities of
$\xi_{11}(zq^{-2s})$ on the sequence $1, q^2,..., q^{2h^\vee-2}$.
Here $\det T_n =4$ (use the first column to compute it from the type $A_n$ case). So again we have no additional relations, and Theorem \ref{t:main} is proved.


\subsection{Type $D_n$}

For $D_n$ we have $r^{\vee} = 1, h^{\vee} = 2(n-1)$. We begin with the case
of even $n$. We have to consider the two half spin representations that we denote $V_{\pm}$.The singularity structure for $\xi_{ij}(z)$ is given by

\begin{align*}
& \xi_{+-}=\xi_{-+} \,: \quad\quad
\begin{matrix}
\bullet^{1} & \bullet^{q^2} & \bullet^{q^4} & \dots & \bullet^{q^{h^{\vee}-2}}&  \bullet^{q^{h^{\vee}}} &  \bullet^{q^{h^{\vee}+2}}& \dots & \bullet^{q^{2h^{\vee}-4}} & \bullet^{q^{2h^{\vee}-2}}\\ 
\hspace{-.25cm}0 & \hspace{-.25cm}1 & \hspace{-.25cm}-1 & \dots & \hspace{-1cm}-1 & \hspace{-.75cm}0 & \hspace{-1cm}-1 & \dots & \hspace{-1cm}-1 & \hspace{-1cm}1
\end{matrix}\\ \\
& \xi_{++}=\xi_{--} \,:\quad\quad
\begin{matrix}
\bullet^{1} & \bullet^{q^2} & \bullet^{q^4} & \dots & \bullet^{q^{h^{\vee}-2}}&  \bullet^{q^{h^{\vee}}} &  \bullet^{q^{h^{\vee}+2}}& \dots & \bullet^{q^{2h^{\vee}-4}} & \bullet^{q^{2h^{\vee}-2}}\\ 
\hspace{-.25cm}2 & \hspace{-.25cm}-1 & \hspace{-.25cm}1 & \dots & \hspace{-1cm}1 & \hspace{-.75cm}-2 & \hspace{-1cm}1 & \dots & \hspace{-1cm}1 & \hspace{-1.25cm}-1
\end{matrix}\\ \\
\end{align*} \vspace{-1.5cm}

\noindent The dots here mean the same thing they meant in the $B_n$ case.

Recall that we have defined indecomposable representations 
$Y(z_1,\dots,z_m|w_1,\dots, w_l)$. The elliptic central character of  
$Y=Y(z_1,\dots,z_m|w_1,\dots,w_l)$ is determined by the pair
\begin{equation}\label{e:ch}
(\xi_{V_+}(z)|_Y,\xi_{V_-}(z)|_Y)=\Big(\prod_{i,j=1}^{m,l} \xi_{++}\big(\frac{z}{z_i}\big)\xi_{+-}\big(\frac{z}{w_j}\big), 
\prod_{i,j=1}^{m,l} \xi_{-+}\big(\frac{z}{z_i}\big)\xi_{--}\big(\frac{z}{w_j}\big)\Big) 
\quad\in\quad \blb C(E)\times \blb C(E)
\end{equation}

\begin{lem}\label{l:linked3}\hfill\vspace{-.25cm}
\begin{enumerate}
\item For any $u\in \blb C^*$,  $Y(z_1,\dots,z_m|w_1,\dots, w_l)$ is linked to 
$Y(z_1,\dots,z_m,  u,  q^{h^{\vee}}u |w_1,\dots, w_l)$ and to 
$Y(z_1,\dots,z_m|w_1,\dots, w_l,  u,  q^{h^{\vee}}u)$. 
In particular, we conclude that $Y(z_1,\dots,z_m|w_1,\dots, w_l)$ is linked to
$Y(z_1,\dots, z_iq^{2h^{\vee}}, \dots ,z_m|w_1,\dots, w_l)$ and to 
$Y(z_1,\dots ,z_m|w_1, \dots, w_jq^{2h^{\vee}}, \dots, w_l)$.

\item $Y(z_1,\dots,z_m,u,uq^{2}|w_1,\dots,w_l,uq^{h^{\vee}},uq^{2+h^{\vee}})$ 
is linked to $Y(z_1,\dots,z_m|w_1,\dots,w_l)$.
\end{enumerate}
\end{lem}

\begin{proof}
The first is proved exactly as in the $B_n$ case since $V_{\pm}^*\cong V_{\pm}(q^{h^{\vee}})$. 
The second follows from the fact that $V_+\otimes V_+(q^{2})$ 
and $V_-\otimes V_-(q^{2})$ have  $V_{n-2}(q)$ 
as a subrepresentation \cite{ChPr91}
and because $V_{n-2}^*\cong V_{n-2}(q^{h^{\vee}})$.
\end{proof}

\obs Adding or removing sequences of the form 
$(u,uq^{2}|uq^{h^{\vee}},uq^{2+h^{\vee}})$, it is easy to show that
$Y(z_1,\dots,z_m,u,uq^{2(2k-1)}|w_1,\dots,w_l,uq^{h^{\vee}},uq^{2(2k-1)+h^{\vee}})$ 
is linked to $Y(z_1,\dots,z_m|w_1,\dots,w_l)$, for $k=1, \dots, (n-2)/2 = (\frac{h^{\vee}}{2}-1)/2$.

Now we have to show that the multiplicative relations between 
$\xi_{++}(zq^{-2s}), \xi_{--}(zq^{-2s})$ and $\xi_{\pm\mp}(zq^{-2s})$,
for $s=0, \dots, \frac{h^{\vee}}{2}-1$, 
are expressed via the transformations of lemma \ref{l:linked3}. 
Since $\xi_{++}(z) = \xi_{--}(z)$ and $\xi_{\pm\mp}(z)=\xi_{\mp\pm}(z)$, we are left to check the relations between $\xi_{++}(zq^{-2s})$ and $\xi_{+-}(zq^{-2s})$.
Consider the group $\blb Z^{\frac{h^{\vee}}{2}}$ and set
\begin{align}
& \xi_{++}(q^{-2s}z) \mapsto v_s := ((-1)^s,(-1)^{s+1}, \dots, -1, 2, -1, \dots, (-1)^s)\\\notag \\
& \xi_{+-}(q^{-2s}z) \mapsto w_s := ((-1)^{s+1},(-1)^s, \dots, 1, 0, 1, \dots, (-1)^{s+1})
\end{align}
Proceeding as in the $B_n$ case, we can show that the $\xi_{++}(zq^{-2s})$,
are multiplicatively independent (the corresponding matrix $T_n$ is of the same form).
The same is true for the $\xi_{+-}(zq^{-2s})$ (using the same arguments). Then,
keeping \eqref{e:ch} in mind, we are left to find coefficients $a_s, a_s'$ such that 
\begin{equation}\label{e:3}
\sum_{s=0}^{\frac{h^{\vee}}{2}-1} a_s v_s = \sum_{s=0}^{\frac{h^{\vee}}{2}-1} a_s' w_s \quad\quad\text{and}\quad\quad \sum_{s=0}^{\frac{h^{\vee}}{2}-1} a_s' v_s = \sum_{s=0}^{\frac{h^{\vee}}{2}-1} a_s w_s
\end{equation}
But it is easy to see that we must have $a_s = a_s'$ and, consequently, that $\sum_s (-1)^s a_s = 0$. We have the following basis of solutions
\begin{equation}
(a_0, a_1, \dots,a_{l-1}, a_l, a_{l+1}, \dots, a_{\frac{h^{\vee}}{2}-1}) = (1, 0, \dots, 0, (-1)^{l+1}, 0, \dots, 0)
\end{equation}
But these are exactly the transformations in the remark after lemma \ref{l:linked3}.

In the odd case the elliptic central character is determined by
its value on one of the half spin representations, $V_n$. The singularity structure for $\xi_{nn}(z)$ is given by the picture

\begin{align*}
\begin{matrix}
\bullet^{1} & \bullet^{q^2} & \bullet^{q^4} & \dots & \bullet^{q^{h^{\vee}-2}}&  \bullet^{q^{h^{\vee}}} &  \bullet^{q^{h^{\vee}+2}}& \dots & \bullet^{q^{2h^{\vee}-4}} & \bullet^{q^{2h^{\vee}-2}}\\ 
\hspace{-.25cm}2 & \hspace{-.25cm}-1 & \hspace{-.25cm}1 & \dots & \hspace{-1cm}-1 & \hspace{-.75cm}0 & \hspace{-1cm}-1 & \dots & \hspace{-1cm}1 & \hspace{-1.25cm}-1
\end{matrix}
\end{align*}

\begin{lem}\label{l:linked4}
$Y(z_1,\dots,z_m,u,uq^{2},uq^{h^\vee}, uq^{h^\vee+2})$ is linked to 
$Y(z_1,\dots,z_m)$. Hence $Y(z_1,\dots, z_m)$ is linked to 
$Y(z_1,\dots, z_{j-1}, z_jq^{2h^{\vee}}, z_{j+1},\dots, z_m)$.
\end{lem}

\begin{proof}
Exactly as the proof of part (b) of lemma \ref{l:linked3}.
\end{proof}

Now we have to check that the relations between $\xi_{nn}(zq^{2s})$,
for $0 \leq s< h^\vee$, are expressed by the transformations of lemma \ref{l:linked4}.
We consider the group $\blb Z^{h^{\vee}}$ and the vectors 
\begin{alignat}3
& v_k &=&\, ((-1)^k,(-1)^{k+1}, \dots, -1, 2, -1, \dots, (-1)^{k+1})\\
& w_k &=&\, ((-1)^{k},(-1)^{k+1}, \dots, -1, 0, -1, \dots, (-1)^{k+1})
\end{alignat}
in $\blb Z^{h^{\vee}/2}$, for $k=0, \dots,
\frac{h^{\vee}}{2}-1=n-2$, 
and assign to $\xi_{nn}(zq^{-2s})$ the vector $u_s \in \blb
Z^{h^{\vee}}$ given by $u_s:=(v_s,w_s)$, if $s< \frac{h^{\vee}}{2}$, or by $u_s:=(w_{s-\frac{h^{\vee}}{2}},v_{s-\frac{h^{\vee}}{2}})$, if $s\geq\frac{h^{\vee}}{2}$. The corresponding matrix $T_n$ is the $h^{\vee}\times h^{\vee}$-matrix of the form
$$T_n = \begin{pmatrix} \bar{T}_n & \dot{T}_n\\ \dot{T}_n & \bar{T}_n \end{pmatrix}$$
where the rows of $\bar{T}_n$ are the vectors $v_k$ and the rows of $\dot{T}_n$ are the $w_k$. By the discussion of the even case we know that the sets $\{v_k\}$ and $\{w_k\}$
are linearly independent. In fact one can show that the rank of $T_n$ is 
$\frac{h^\vee}{2}+1 = n$. The linear relations between the two sets of vectors 
are similar to that of the even case and have the form,
\begin{equation}
u_0+u_{\frac{h^{\vee}}{2}} +
(-1)^{k+1}(u_k+u_{\frac{h^{\vee}}{2}+k}) =0,
 \quad\quad k=1, \dots, \frac{h^{\vee}}{2}-1
\end{equation}
Similarly to the remark after lemma \ref{l:linked3}, one shows that these relations are all 
obtained from the one in lemma \ref{l:linked4} (here will be used that $V_n^*\cong V_{n-1}(q^{h^{\vee}})$).


\subsection{Type $E_n$}

\noindent $\bullet\,E_6$

The dual Coxeter number is $h^{\vee} = 12$ and $r^{\vee}=1$.
The black node corresponds to one of the 27-dimensional fundamental representations. The singularity structure for $\xi_{11}(z)$ is

$$\begin{matrix}
\bullet^{1} & \bullet^{q^2} & \bullet^{q^4} & \bullet^{q^6} & \bullet^{q^8} & \bullet^{q^{10}} &  \bullet^{q^{12}}\\
\hspace{-.25cm}2 & \hspace{-.25cm}-1 & \hspace{-.35cm}0 & \hspace{-.25cm}1 & \hspace{-.25cm}-1 & \hspace{-.35cm}0 & \hspace{-.35cm}0
\end{matrix}\hspace{.25cm}
\begin{matrix}
\bullet^{q^{14}} & \bullet^{q^{16}} & \bullet^{q^{18}} & \bullet^{q^{20}} & \bullet^{q^{22}}\\  
\hspace{-.35cm}0 & \hspace{-.25cm}-1 & \hspace{-.25cm}1 & \hspace{-.35cm}0 & \hspace{-.25cm}-1
\end{matrix}$$

\begin{lem}\label{l:linkede6}\hfill\vspace{-.25cm}
\begin{enumerate}
\item $Y(z_1,\dots,z_m,w,wq^{8},wq^{16})$ is linked to $Y(z_1,\dots,z_m)$. 
In particular $Y(z_1,\dots, z_m)$ is linked to 
$Y(z_1,\dots, z_{j-1}, z_jq^{24}, z_{j+1},\dots, z_m)$.

\item $Y(z_1,\dots,z_m,w,wq^2,wq^{10},wq^{12})$ is linked to $Y(z_1,\dots,z_m,wq^6)$.
\end{enumerate}
\end{lem}

\begin{proof}
Item (a) follows since  
$V_5(q^4)$ is a subrepresentation of $V_1\otimes V_1(q^8)$ \cite{ChPr91} and 
because $V_1^*\cong V_5(q^{12})$. For item (b) we will use a sequence of subrepresentations, all of them computed in \cite{ChPr91}.
By part (a), $V_1(q^2)\otimes V_1(q^{10})$ has 
$V_5(q^{6})$ as a subrepresentation and, consequently,
$V_1\otimes V_1(q^2)\otimes V_1(q^{10})\otimes V_1(q^{12})$ has 
$V_1\otimes V_5(q^6)\otimes V_1(q^{12})$ as a subrepresentation. 
Now, $V_6(q^{3})$ is a subrepresentation of $V_1\otimes V_5(q^{6})$ and we get that
$V_1\otimes V_6(q^9)$ is a subrepresentation of $V_1\otimes V_5(q^6)\otimes V_1(q^{12})$.
Finally we get $V_1(q^6)$ as subrepresentation of $V_1\otimes V_6(q^9)$.
\end{proof}

\obs All these subrepresentation relations can be obtained studying tensor product 
decompositions and the singularity structure of the elliptic central character. 
This will be the procedure to prove theorem \ref{t:generator} 
for the remaining $\lie g$.

The combinatorics part is in the group $\blb Z^{12}$. We consider the vectors $v_0, \dots, v_{11}$ obtained from
$$v_0 = (2,-1,0,1,-1,0,0,0,-1,1,0,-1)$$
by cyclic permutation of the coordinates, corresponding to $\xi_{11}(zq^{-2s})$, where the entries of $v_s$ have the usual meaning.
The rank of the {\small(12x12)}-matrix $T_n$ thus obtained is 6. We first eliminate the linear 
relations related to item (a) of lemma \ref{l:linkede6}.
They are
\begin{equation}
v_k+v_{k+4}+v_{k+8}=0
\end{equation}
for $k=0,1,2,3$.
Removing $v_8,\dots, v_{11}$ from our set of vectors, we check that the remaining vectors 
$v_0, \dots, v_7$ satisfy the following linear relations
\begin{equation}
v_{k+3}=v_k+v_{k+1}+v_{k+5}+v_{k+6}
\end{equation}
for $k=0,1$. This is exactly part (b) of lemma \ref{l:linkede6}. Since the vectors 
$v_0,\dots,v_5$ are linearly independent, the proof is complete.

\obs The rank and the kernel of the corresponding matrix $T_n$ for $\lie e_n$ and $\lie f_4$
were computed using the computer software {\it Mathematica}.

\vspace{.5cm}\noindent $\bullet\,E_7$

Here $h^{\vee} = 18$. The black node corresponds to the
affinization of the 56-dimensional fundamental representation. We
list the singularities of $\xi_{11}$ and $\xi_{16}$. These are the two we need to prove Theorem \ref{t:main}. To prove Theorem \ref{t:generator} we will need to analyse other $\xi_{1j}(z)$, but they can be read from the matrices $M(q)$ in the appendix.

\begin{align*}
& \xi_{11} \,: \quad 
\begin{matrix}
\bullet^{1} & \bullet^{q^2} & \hspace{1cm}\bullet^{q^8} & \bullet^{q^{10}} & \hspace{1cm} \bullet^{q^{16}} & \bullet^{q^{18}} & \bullet^{q^{20}} & \hspace{1cm}\bullet^{q^{26}} & \bullet^{q^{28}} & \hspace{1cm}\bullet^{q^{34}}\\  
\hspace{-.25cm}2 & \hspace{-.25cm}-1 & \hspace{1cm}\hspace{-.25cm}1 & \hspace{-.25cm}-1 & \hspace{1cm}\hspace{-.35cm}1 & \hspace{-.35cm}-2 & \hspace{-.35cm}1 & \hspace{1cm}\hspace{-.25cm}-1 & \hspace{-.25cm}1 & \hspace{1cm}\hspace{-.25cm}-1
\end{matrix}\\ \\
& \xi_{16} \,: \qquad\qquad 
\begin{matrix}
\bullet^{q^5} & \bullet^{q^7} & \hspace{1cm}\bullet^{q^{11}} & \bullet^{q^{13}} & \hspace{2cm}\bullet^{q^{23}} & \bullet^{q^{25}} & \hspace{1cm}\bullet^{q^{29}} & \bullet^{q^{31}}\\
\hspace{-.25cm}1 & \hspace{-.25cm}-1 & \hspace{1cm}\hspace{-.25cm}1 & \hspace{-.25cm}-1 & \hspace{2cm}\hspace{-.25cm}-1 & \hspace{-.25cm}1 & \hspace{1cm}\hspace{-.25cm}-1 & \hspace{-.25cm}1 
\end{matrix}
\end{align*}

We begin by proving Theorem \ref{t:generator}. We denote the
finite dimensional representation of $U_q(\lie e_7)$ with highest
weight $\lambda$ by $V_{\lambda}^f$. It is known that $V_1^f$ is
affinizable, i.e., the $U_q(\hlie e_7)$ fundamental
representation $V_1$ is isomorphic to $V_1^f$ as a $U_q(\lie e_7)$-module. Now

\begin{equation}\label{e:V1V1}
V_1^f\otimes V_1^f \cong V_{2\omega_1}^f\oplus V_2^f \oplus V_6^f\oplus \blb C
\end{equation}
Using Corollary \ref{c:poles} we see that the normalized R-matrix
is not invertible at $q^{-2}, q^{-10}$ and $q^{-18}$. The
associated elliptic central characters for these cases are,
respectively, the same as the 
elliptic central characters as $V_2(q), V_6(q^5)$ and $\blb C$, respectively. In fact these are subrepresentations : by \eqref{e:V1V1} and the associated elliptic central character, the kernel of $P\bar R_{V_1,V_1}(q^{-2})$ must be of the form $V_2(q^{1+36k})$ for some $k \in \blb Z$. To conclude that $k=0$, one use corollary \ref{c:smallshift} with $i=j=1, r=2$ and
$m=1$. The story for $V_6(q^5)$ is exactly the same. The affinization of $V_6^f$ (which corresponds to the adjoint representation of $\lie e_7$) is isomorphic to $V_6^f\oplus\blb C$ and thus we have
$$V_6 \otimes V_6(w) \cong V_{2\omega_6}^f \oplus V_5^f \oplus V_2^f \oplus 3V_6^f \oplus 2\blb C$$ 
The R-matrix is non invertible when $w=q^2, q^8, q^{12}, q^{18}$. We proceed as before to check that $V_5(q)$ is a subrepresentation of $V_6\otimes V_6(q^2)$. This is the only new fundamental representation we get here. Then we go to
$$V_1\otimes V_6(w) \cong V_{\omega_1+\omega_6}^f\oplus V_7^f\oplus 2V_1$$ 
and the R-matrix is non invertible at $w = q^7, q^{13}$. The new fundamental representation we get is $V_7(q^4)$ as a subrepresentation of $V_1\otimes V_6(q^7)$. 

\obs For the next tensor product we will need the following important observations. All the possible affinizations of $V_{\omega_1+\omega_6}^f$ (denoted by 
$V_{\omega_1+\omega_6}[w](u)$) occur as (shifts) of subquotients of $V_1(u)\otimes V_6(wu)$. Therefore, all possible $\xi_{V_1}(z)|_{V_{\omega_1+\omega_6}[w](u)}$ are of the form $\xi_{11}(z/u)\xi_{16}(z/wu)$.
                
We now study the tensor product $V_1\otimes V_2(w)$. The affinization of $V_2^f$ is isomorphic to $V_2^f\oplus V_6^f\oplus \blb C$. We have
$$V_1\otimes V_2(w) \cong V_{\omega_1+\omega_2}^f \oplus V_3^f \oplus 2V_{\omega_1+\omega_6}^f \oplus 2V_7^f \oplus 3V_1^f$$
The values of $w$ where $V_1\otimes V_2(w)$ has a subrepresentation not containing the highest weight component are $q^3, q^{11}$ and $q^{17}$ and, at these points, the corresponding elliptic central characters are those of $V_3(q^2), V_7(q^8)$ and $V_1(q^{16})$. By the last remark we know that all possible elliptic central characters of $V_{\omega_1+\omega_6}[w](u)$ at $V_1$ are of the form $\xi_{11}(z/u)\xi_{16}(z/wu)$. One checks now that this will never produce the elliptic central character of $V_3(x)$  for any $x$ and, therefore, we do get $V_3(q^2)$ as a subrepresentation. To complete the proof of Theorem \ref{t:generator} we need to get $V_4$. This is done as before using that
$$V_1^f\otimes V_3^f \cong V_{\omega_1+\omega_3}^f \oplus
V_4^f\oplus V_{\omega_2+\omega_6}^f\oplus V_{\omega_1+\omega_7}^f
\oplus V_5^f \oplus V_2^f$$
Thus
 $V_4(q^3)$ occurs as subrepresentation of $V_1\otimes V_3(q^4)$.

\obs The tensor product decompositions for $\lie e_7$ and $\lie e_8$ were computed using the
computer package {\it LiE} (http://wwwmathlabo.univ-poitiers.fr/~maavl/LiE).

Now we proceed with the proof of Theorem \ref{t:main}.

\begin{lem}\label{l:linkede7}\hfill\vspace{-.25cm}
\begin{enumerate}
\item For any $w\in \blb C^*$, $Y(z_1,\dots, z_m)$ is linked to 
$Y(z_1,\dots, z_m, w, wq^{18})$. In particular $Y(z_1,\dots, z_m)$ is linked to $Y(z_1,\dots, z_{j-1}, z_jq^{36}, z_{j+1},\dots, z_m)$.

\item $Y(z_1,\dots,z_m,w,wq^2,wq^{12},wq^{14},wq^{24},wq^{26})$ is linked to   
$Y(z_1,\dots,z_m)$.
\end{enumerate}
\end{lem}

\begin{proof}
Part (a) is clear from $V_1^*\cong V_1(q^{18})$. Let us prove
(b). Using that $V_6(q^5)$ is a subrepresentation of $V_1\otimes V_1(q^{10})$ we find that $V_1 \otimes V_6(q^7)\otimes V_6(q^{19})$ is a subrepresentation of $V_1\otimes V_1(q^2)\otimes V_1(q^{12})\otimes V_1(q^{14})\otimes V_1(q^{24})$. Now use that $V_6(q^6)$ is a subrepresentation of $V_6\otimes V_6(q^{12})$ to get that $V_1\otimes V_6(q^{13})$ is a subrepresentation of
$V_1 \otimes V_6(q^7)\otimes V_6(q^{19})$. Since $V_1(q^8)$ is a
subrepresentation of $V_1\otimes V_6(q^{13})$ and
$V_1(q^{26})\cong V_1(q^8)^*$, we have $\Bbb C\subset 
V_1\otimes V_1(q^2)\otimes V_1(q^{12})\otimes V_1(q^{14})\otimes
V_1(q^{24})\otimes V_1(q^{26})$, so we are done.
\end{proof}

We now consider the group $\blb Z^9$ and the vectors $v_0, \dots, v_8$ corresponding to $\xi_{11}(zq^{-2s})$ as usual
(the entries are the orders of the singularities of
$\xi_{11}(zq^{-2s})$ on the sequence $1, q^2,..., q^{16}$).
The matrix $T_n$ has rank 7 and the non-trivial linear relations are
\begin{equation}\label{e:e7}
v_k+v_{k+1}+v_{k+6}+v_{k+7} - v_{k+3}-v_{k+4}=0
\end{equation}
for $k=0,1$. These relations are implemented by part (b) of
lemma \ref{l:linkede7}, 
so Theorem \ref{t:main} is proved. 

\vspace{.5cm}\noindent $\bullet\,E_8$

For $\lie e_8$ $h^{\vee} = 30$ and $r^{\vee}=1$. The black node corresponds to the affinization of the adjoint representation. The singularity structures for $\xi_{11}$ and $\xi_{17}$ are respectively

\begin{align*}
& \begin{matrix}
\bullet^{1} & \bullet^{q^2} & \hspace{.5cm}\bullet^{q^{10}} & \bullet^{q^{12}} & \hspace{.5cm} \bullet^{q^{18}} & \bullet^{q^{20}} & \hspace{.5cm}\bullet^{q^{28}} & \bullet^{q^{30}} & \bullet^{q^{32}}\\  
\hspace{-.25cm}2 & \hspace{-.25cm}-1 & \hspace{.5cm}\hspace{-.25cm}1 & \hspace{-.25cm}-1 & \hspace{.5cm}\hspace{-.35cm}1 & \hspace{-.35cm}-1 & \hspace{.5cm}\hspace{-.35cm}1 & \hspace{-.25cm}-2 & \hspace{-.25cm}1
\end{matrix}\hspace{.75cm}
\begin{matrix}
\bullet^{q^{40}} & \bullet^{q^{42}} & \hspace{.5cm} \bullet^{q^{48}} & \bullet^{q^{50}} & \hspace{.5cm}\bullet^{q^{58}}\\  
\hspace{-.25cm}-1 & \hspace{-.25cm}1 & \hspace{.5cm}\hspace{-.35cm}-1 & \hspace{-.35cm}1 & \hspace{.5cm}\hspace{-.35cm}-1
\end{matrix}\\ \\
& \quad\begin{matrix}
\bullet^{q^6} & \bullet^{q^8} & \bullet^{q^{12}} & \bullet^{q^{14}} & \bullet^{q^{16}} & \bullet^{q^{18}} & \bullet^{q^{22}} & \bullet^{q^{24}}\\  
\hspace{-.25cm}1 & \hspace{-.25cm}-1 & \hspace{-.25cm}1 & \hspace{-.25cm}-1 & \hspace{-.35cm}1 & \hspace{-.35cm}-1 & \hspace{-.35cm}1 & \hspace{-.25cm}-1
\end{matrix}\hspace{.5cm}
\begin{matrix}
\bullet^{q^{36}} & \bullet^{q^{38}} & \bullet^{q^{42}} & \bullet^{q^{44}} & \bullet^{q^{46}} & \bullet^{q^{48}} & \bullet^{q^{52}}& \bullet^{q^{54}}\\  
\hspace{-.25cm}-1 & \hspace{-.25cm}1 & \hspace{-.35cm}-1 & \hspace{-.35cm}1 & \hspace{-.35cm}-1 & \hspace{-.35cm}1 & \hspace{-.35cm}-1 & \hspace{-.35cm}1
\end{matrix}\\ 
\end{align*}

Theorem \ref{t:generator} is proved as in the $\lie e_7$ case using the following relations

\begin{gather}\notag
V_2(q) \subset V_1\otimes V_1(q^2)  \quad\quad V_7(q^{6}) \subset V_1\otimes V_1(q^{12}) \quad\quad V_1(q^{10}) \subset V_1\otimes V_1(q^{20})\\\label{e:e8}
V_8(q^{5}) \subset V_1\otimes V_7(q^{8}) \quad\quad V_6(q) \subset V_7\otimes V_7(q^2) \quad\quad V_3(q^2)  \subset V_1\otimes V_2(q^3)\\\notag
V_4(q^{3}) \subset V_1\otimes V_3(q^{4})  \quad\quad   V_5(q^{4}) \subset V_1\otimes V_4(q^{5})
\end{gather}

The tensor product decomposition for fundamental representations of $\lie e_8$ involves much
longer expressions than those we had to deal in the $\lie e_7$ case, so we will not write 
them down completely. We will just write two steps, since we will need them to prove 
Theorem \ref{t:main}. We sketch the remaining steps in the appendix. 
The first step is the tensor product $V_1\otimes V_1(w)$. 
Since $V_1 \cong V_1^f\oplus \blb C$ we have
$$V_1\otimes V_1(w) \cong V_{2\omega_1}^f \oplus V_2^f\oplus V_7^f\oplus 3V_1^f\oplus 2\blb C$$
From this we prove the first line of \eqref{e:e8} and check that
the elliptic central character (at $V_1$) of the affinizations
$V_{2\omega_1}[w](u)$ coincides with that of a fundamental
representation only when $w \in \{q^{\pm 2}, q^{\pm 12}, q^{\pm
  20}\}$ (modulo $q^{60}$), when it has the elliptic central
character of $V_2(q^{\pm 1}), V_7(q^{\pm 6})$ or $V_1(q^{\pm
  10})$, respectively. The proof that these are really
subrepresentations when $w \in \{q^{2}, q^{12}, q^{20}\}$, is
analogous to the one we did for $\lie e_7$.

In the second step we consider $V_1\otimes V_7(w)$. Since $V_7 \cong V_7^f\oplus V_1^f\oplus \blb C$ we have
$$V_1\otimes V_7(w) \cong V_{\omega_1+\omega_7}^f \oplus V_8^f\oplus V_{2\omega_1}^f\oplus  2V_2^f\oplus 3V_7^f\oplus 4V_1^f\oplus 2\blb C$$
The normalized R-matrix is not invertible at $w=q^8, q^{14}, q^{18}, q^{24}$ and the corresponding elliptic central characters coincide with the ones of $V_8(q^5), V_2(q^9), V_7(q^{12})$ and $V_1(q^{18})$, respectively. Since all other components do not have the elliptic central character of $V_8$, we conclude that we really get an affinization of $V_8^f$ as subrepresentation, and it must be $V_8(q^5)$. We will come back to this tensor product in the proof of lemma \ref{l:linkede8} below.
                                        
\begin{lem}\label{l:linkede8}\hfill\vspace{-.25cm}
\begin{enumerate}
\item For any $w\in \blb C^*$, $Y(z_1,\dots, z_m)$ is linked to 
$Y(z_1,\dots, z_m, w, wq^{30})$. In particular $Y(z_1,\dots, z_m)$ is linked to $Y(z_1,\dots, z_{j-1}, z_jq^{60}, z_{j+1},\dots, z_m)$.

\item $Y(z_1,\dots, z_m,w,wq^{20},wq^{40})$ is linked to $Y(z_1,\dots, z_m)$.

\item $Y(z_1,\dots, z_m,w,wq^{12},wq^{24},wq^{36},wq^{48})$ is linked to $Y(z_1,\dots, z_m)$.
\end{enumerate}
\end{lem}

\begin{proof}
Part (a) is clear from $V_1^*\cong V_1(q^{30})$. Part (b) follows since $V_1(q^{10})$ 
is a subrepresentation of $V_1\otimes V_1(q^{20})$. 
Part (c) is more delicate. First $V_7(q^6) \subset V_1\otimes V_1(q^{12})$ and,
consequently,  $V_1\otimes V_7(q^{18})\otimes V_7(q^{42})$ is a subrepresentation of 
$V_1\otimes V_1(q^{12})\otimes V_1(q^{24})\otimes V_1(q^{36})\otimes V_1(q^{48})$. 
Let us introduce the notation $Y(z_1,\dots,z_m|u_1,\dots,u_l)$ meaning a non-resonant reordering of 
$V_1(z_1)\otimes \dots\otimes V_1(z_m)\otimes V_7(u_1)\otimes \dots\otimes V_7(u_l)$. 
Then we have shown that
 $Y(z_1,\dots z_m,w|wq^{18},wq^{42})$ is a subrepresentation of 
$Y(z_1,\dots, z_m,w,wq^{12},wq^{24},wq^{36},wq^{48})$

Go back to the tensor product $V_1\otimes V_7(w)$. We cannot
conclude whether we have $V_7(q^{12})$ as a subrepresentation of
$V_1\otimes V_7(q^{18})$. If it was true the proof would be
completed, as $V_7^*\cong V_7(q^{30})$. 
If it is not true, then we must have, as subrepresentation, an irreducible affinization of $V_{2\omega_1}^f$ with the elliptic central character of 
$V_7(q^{12})$.
From the analysis of $V_1\otimes V_1(w)$ we know that  this affinization must be a quotient 
of the form
$$V_1(q^{6})\otimes V_1(q^{18+60j}) \twoheadrightarrow
V_{2\omega_1}[q^{12+60j}](q^6)$$
Let $Y(\overrightarrow{z},(w|wq^{18}),wq^{42})$, where 
$\overrightarrow{z}=(z_1,\dots, z_m)$, denote the corresponding subrepresentation 
of the original $Y(\overrightarrow{z},w,wq^{12},wq^{24},wq^{36},wq^{48})$. We have shown
$$Y(\overrightarrow{z},wq^{6},wq^{18+60j}|wq^{42}) \twoheadrightarrow
Y(\overrightarrow{z},(w|wq^{18}),wq^{42}) \hookrightarrow
Y(\overrightarrow{z},w,wq^{12},wq^{24},wq^{36},wq^{48})$$
We complete the proof with the following diagram
$$Y(\overrightarrow{z}) \hookrightarrow Y(\overrightarrow{z}|wq^{12},wq^{42}) 
\hookrightarrow Y(\overrightarrow{z},wq^6,wq^{18}|wq^{42}) \sim 
Y(\overrightarrow{z},wq^{6},wq^{18+60j}|wq^{42})$$
The first inclusion follows from $V_7^*\cong V_7(q^{30})$ and the second since 
$V_7(q^{12}) \hookrightarrow V_1(q^6)\otimes V_1(q^{18})$. The symbol $\sim$ denotes a
linking relation and is immediate from part (a).
\end{proof}

As usual, now we define vectors $v_0, \dots, v_{14}$ in $\blb Z^{15}$, encoding the orders of the singularities of
$\xi_{11}(zq^{-2s})$ on the sequence $1, q^2,..., q^{28}$.
The corresponding matrix $T_n$ has rank 8. We find the following (linearly independent) relations
\begin{align}\label{e1}
v_j+v_{j+10}-v_{j+5} = 0 & \quad\quad\text{for }j=0,1,2,3,4\\\label{e2}
v_k+v_{k+6}+v_{k+12} - v_{k+3}-v_{k+9} = 0 & \quad\quad\text{for }k=0,1
\end{align}
Relation \eqref{e1} is part (b) of lemma \ref{l:linkede8},
while \eqref{e2} is part (c).


\subsection{$F_4$}

For $\lie f_4$ we have $r^{\vee} = 2$ and $h^{\vee} = 9$. The black node corresponds to $V_1$, the affinization of the 26-dimensional representation of $\lie f_4$. We list the singularity structure only for $\xi_{11}(z)$.

$$\begin{matrix}
\bullet^{1} & \bullet^{q^2}& \bullet^{q^4} & \bullet^{q^6} & \bullet^{q^{8}} & \bullet^{q^{10}} & \bullet^{q^{12}} & \bullet^{q^{14}} & \bullet^{q^{16}}\\  
\hspace{-.25cm}2 & \hspace{-.25cm}-1 & \hspace{-.35cm}0 & \hspace{-.25cm}1 & \hspace{-.25cm}-1 & \hspace{-.25cm}1 & \hspace{-.25cm}-1 & \hspace{-.35cm}0 & \hspace{-.35cm}1
\end{matrix}\,
\begin{matrix}
\bullet^{q^{18}} & \bullet^{q^{20}}& \bullet^{q^{22}} & \bullet^{q^{24}} & \bullet^{q^{26}} & \bullet^{q^{28}} & \bullet^{q^{30}} & \bullet^{q^{32}} & \bullet^{q^{34}}\\  
\hspace{-.25cm}-2 & \hspace{-.25cm}1 & \hspace{-.35cm}0 & \hspace{-.25cm}-1 & \hspace{-.25cm}1 & \hspace{-.25cm}-1 & \hspace{-.25cm}1 & \hspace{-.35cm}0 & \hspace{-.25cm}-1
\end{matrix}$$

$V_1^f$ is affinizable and we have

$$V_1 \otimes V_1(w) \cong V_{2\omega_1}^{f} \oplus V_2^{f} \oplus V_4^{f} \oplus V_1^{f} \oplus \blb C$$
By Corollary \ref{c:poles}, $\bar{R}_{V_1,V_1}(z)$ is not invertible at
$z=q^{-2}, q^{-8}, q^{-12}, q^{-18}$. Arguing as in the $\lie e_7$ case we prove 
that the subrepresentations of $V_1\otimes V_1(w)$, for $w=q^{2}, q^{8}, q^{12}, q^{18}$
are, respectively,  $V_2(q), V_4(q^{4}), V_1(q^{6})$ and \blb C. The affinization of the adjoint representation $V_4$ is isomorphic to $V_4^f\oplus \blb C$ and we have
$$V_4\otimes V_4(w) \cong V_{2\omega_4}^f \oplus V_3^f\oplus V_{2\omega_1}^f\oplus 3V_4^f \oplus 2\blb C$$
As before we get that $V_3(q^2)$ is a subrepresentation of $V_4\otimes V_4(q^4)$ and Theorem \ref{t:generator} is proved.

To prove Theorem \ref{t:main} we state

\begin{lem}\label{l:linkedf4}\hfill\vspace{-.25cm}
\begin{enumerate}
\item For any $w\in \blb C^*$, $Y(z_1,\dots, z_m)$ is linked to 
$Y(z_1,\dots, z_m, w, wq^{18})$. In particular $Y(z_1,\dots, z_m)$ is linked to $Y(z_1,\dots, z_{j-1}, z_jq^{36}, z_{j+1},\dots, z_m)$.

\item $Y(z_1,\dots, z_m,w,wq^{12},wq^{24})$ is linked to $Y(z_1,\dots, z_m)$.
\end{enumerate}
\end{lem} 

\begin{proof}
Part (a) follows from $V_1^*\cong V_1(q^{18})$. Since $V_1(q^6)$ is a subrepresentation of $V_1\otimes V_1(q^{12})$, part (b) is proved.
\end{proof}

We define $v_0,\dots, v_8 \in \blb Z^9$ corresponding to $\xi_{11}(zq^{-2s})$. The rank of $T_n$ is 6. The linear relations give exactly part (b) of lemma \ref{l:linkedf4}.
\begin{equation}
v_k+v_{k+6} - v_{k+3} = 0 \quad\quad\text{for } k=0,1,2 
\end{equation}


\subsection{$G_2$}

In this case $r^{\vee} = 3$ and $h^{\vee} = 4$. The zeros and poles are given by
\begin{align*}
& \xi_{11}(z)\quad : \qquad\qquad
\begin{matrix}
\bullet^{1} & \bullet^{q^2}& \bullet^{q^4} & \bullet^{q^8} & \bullet^{q^{10}} & \bullet^{q^{12}} & \bullet^{q^{14}} & \bullet^{q^{16}} & \bullet^{q^{20}}& \bullet^{q^{22}}\\  
\hspace{-.25cm}2 & \hspace{-.25cm}-1 & \hspace{-.25cm}1 & \hspace{-.25cm}-1 & \hspace{-.25cm}1 & \hspace{-.25cm}-2 & \hspace{-.25cm}1 & \hspace{-.25cm}-1 & \hspace{-.25cm}1 & \hspace{-.25cm}-1
\end{matrix}\\ \\
& \xi_{12}(z)\quad : \qquad\qquad\qquad
\begin{matrix}
\bullet^{q} &  \bullet^{q^5}& \bullet^{q^7} & \bullet^{q^{11}} & \bullet^{q^{13}} & \bullet^{q^{17}} & \bullet^{q^{19}} & \bullet^{q^{23}}\\  
\hspace{-.25cm}1 & \hspace{-.25cm}1 & \hspace{-.25cm}-1 & \hspace{-.25cm}-1 & \hspace{-.25cm}-1 & \hspace{-.25cm}-1 & \hspace{-.25cm}1 & \hspace{-.25cm}1
\end{matrix}
\end{align*} 
The black node corresponds to $V_1$, the affinization of the 7-dimensional 
representation of $\lie g_2$. $V_1$ is affinizable and
$$V_1\otimes V_1(w) \cong V_{2\omega_1}^f\oplus V_2^f\oplus V_1^f\oplus \blb C$$
Subrepresentations are obtained at $w=q^2, q^8, q^{12}$ and they are $V_2(q), V_1(q^4)$ and \blb C, respectively. This proves Theorem \ref{t:generator}.

The usual lemma is immediate.

\begin{lem}\label{l:linkedg2}\hfill\vspace{-.25cm}
\begin{enumerate}
\item For any $w\in \blb C^*$, $Y(z_1,\dots, z_m)$ is linked to 
$Y(z_1,\dots, z_m, w, wq^{12})$. In particular $Y(z_1,\dots, z_m)$ is linked to $Y(z_1,\dots, z_{j-1}, z_jq^{24}, z_{j+1},\dots, z_m)$.

\item $Y(z_1,\dots, z_m,w,wq^{8},wq^{16})$ is linked to $Y(z_1,\dots, z_m)$.
\end{enumerate}
\end{lem} 

Then we consider the vectors
 encoding the sigularity structures of $\xi_{11}(zq^{-2s})$ in the sequence $1,q^2,\dots, q^{10}$ :
\begin{gather*}
v_0 = (2,-1,1,0,-1,1) \quad v_1=(-1,2,-1,1,0,-1)  \quad v_2=(1,-1,2,-1,1,0)\\
v_3 = (0,1,-1,2,-1,1) \quad v_4=(-1,0,1,-1,2,-1)  \quad v_5=(1,-1,0,1,-1,2)
\end{gather*}
and observe that $v_{2+k}=v_k+v_{4+k}$, for $k=0,1$, which is part (b) of the lemma. One easily checks that $v_0,\dots, v_3$ are linearly independent to complete the proof.

\appendix
\section*{Appendix}

\section{Sketch of the Proof of Theorem \ref{t:generator} for $E_8$} 

The idea of the proof is the following (recall the definition of the set $\cal P_{ij}$ in corollary \ref{c:poles}).

\noindent{\bf 1.} Compute the $U_q(\lie g)$-tensor product decomposition of $V_i\otimes V_j(w)$, 
for already obtained $V_i,V_j$, and check that a new $V_k^f$ occurs.

\noindent{\bf 2.} Compute all possible elliptic central characters (ECC) of 
$V_i\otimes V_j(w)$ at $V_1$ and check that it will coincide with the 
ECC of the new $V_k$ for some $w \in P_{ij}$.

\noindent{\bf 3.} Check that all affinizations of the remaining components never have
the ECC of $V_k$ to conclude that $V_k(wq^s)$ is indeed a subrepresentation of
$V_i\otimes V_j(w)$, for some $s\in \blb Z$.

It becomes clear that we need the list of the singularity structures for $\xi_{ij}$. 
They can be obtained directly from the matrix $M(q)$ given in the section below. 
The following lemma is helpful. 

\begin{lem}\label{l:fECC}
$V_i\otimes V_j(w)$ has ECC of some $V_k(u)$ if and only if $w$ belongs 
(modulo $q^{60})$ to $\cal P_{ij}^{\pm 1}$.
\end{lem}

\begin{proof}
First observe that the singularities of $\xi_{ij}$ have order $\pm 1$ or $\pm 2$, and that 
if order $\pm 2$ occurs, it occurs only once and, in that case, $i=j$.
Hence, for $\xi_{ii}(z)\xi_{ij}(z/w)$, coincide with some $\xi_{ik}(zu)$, 
the zero of order $2$ in $\xi_{ii}(z)$ must be combined with a pole of $\xi_{ij}(z/w)$.
\end{proof}

Using the periodicity properties of $\xi_{ij}$, we see that to find all possible 
``fundamental'' ECC of $V_i\otimes V_j(w)$, we just need to take $w\in P_{ij}$. In Table 2 we give this list for the pairs $(i,j)$ that we will need

\vspace{.5cm}{\centerline{\bf Table 2}}

\hspace*{-.65cm}\begin{tabular}{|m{.75cm}|m{3.1cm}|m{2.5cm}|m{.75cm}|m{3.75cm}|m{3.1cm}|}
\hline Pair &   $\mathcal{P}_{ij}$                     & ECC of & Pair & $\mathcal {P}_{ij}$              & ECC of         \\
\hline (1,1)       & {\small $\{q^{2},q^{12},q^{20},q^{30}\}$}     & $V_2,{V}_7,{V}_1,\blb C$ 
& (2,2) & {\small $\{q^{4},q^{12},q^{14},q^{20},q^{22},q^{30}\}$} & ${V}_4,{V}_6,{V}_3,{V}_2,{V}_7,\blb C$ \\
\hline (1,2)       & {\small $\{q^{3},q^{13},q^{21},q^{29}\}$}     & ${V}_3,{V}_8,{V}_7,{V}_1$
& (2,7)             & {\small$\{q^{9},q^{19},q^{25}\}$}                  & ${V}_6,{V}_2,{V}_1$ \\
\hline (1,3)       & {\small $\{q^{4},q^{14},q^{22},q^{28}\}$}     & ${V}_4,{V}_6,{V}_8,{V}_2$
& (2,8)             & {\small $\{q^{16},q^{25}\}$}                            & ${V}_3,{V}_1$ \\
\hline (1,4)       & {\small $\{q^{5},q^{23},q^{27}\}$}                 & ${V}_5,{V}_6,{V}_3$
& (3,7)             & {\small $\{q^{16},q^{25}\}$}                          & ${V}_4,{V}_6,{V}_8,{V}_2$\\
\hline (1,6)       & {\small $\{q^{9},q^{19},q^{25}\}$}                  & $V_4,V_3,V_8$
& (7,7) & {\small $\{q^{2},q^{8},q^{14},q^{20},q^{24},q^{30}\}$} & $V_6,V_8,V_2,V_7,V_1,\blb C$ \\
\hline (1,7) & {\small $\{q^{8},q^{14},q^{18},q^{24}\}$}    & $V_8,V_2,V_7,V_1$
& (7,8)             & {\small$\{q^{5},q^{23},q^{27}\}$}                   & $V_4,V_7,V_1$\\
\hline (1,8) & {\small $\{q^{7},q^{11},q^{17},q^{21},q^{25}\}$} & $V_6,V_3,V_8,V_2,V_7$ & & & \\
\hline
\end{tabular}

\vspace{.5cm}
Recall that we began studying $V_1\otimes V_1(w)$ and obtained 
$V_2$ and $V_7$ as subrepresentations. Then we considered $V_1\otimes V_7(w)$ 
and obtained $V_8$ as subrepresentation. It is clear that multiplicities do not affect our arguments, 
so we will list the components and will add a sign ``$\oplus \rm{mult}$'' at the end to indicate that
the remaining factors have already been listed. We continue with

$$V_7\otimes V_7(w) \cong V^f_{2\omega_7}\oplus V^f_6 \oplus V^f_{\omega_1+\omega_7}
\oplus V^f_8\oplus V^f_{2\omega_1} \oplus V^f_2\oplus V^f_7\oplus V^f_1\oplus \blb C \oplus \rm{mult}$$
$$V_1\otimes V_2(w) \cong V^f_{\omega_1+\omega_2}\oplus V^f_3 \oplus V^f_{\omega_1+\omega_7} \oplus V^f_8\oplus V^f_{2\omega_1} \oplus V^f_2\oplus V^f_7\oplus V^f_1\oplus
 \blb C \oplus \rm{mult}$$
Looking at Table 2 we see that affinizations of $V_{\omega_1+\omega_7}$ and 
$V_{2\omega_1}$ never have ECC of $V_6$ nor of
$V_3$
(up to shift). 
Therefore the kernel of 
$\bar{R}_{V_7,V_7}(q^2)$ must be an affinization of $V^f_6$, while the kernel of
$\bar{R}_{V_1,V_2}(q^3)$ must be an affinization of $V^f_3$. We will not use the knowledge of
the ``$U_q(\lie g)$-tail'' of $V_3$ (i.e. the
representations of $U_q(\lie g)$ 
added to $V^f_3$ to obtain $V_3$). So, for the next tensor product, we will assume
that $V_3$ has the longest possible ``tail''. Then we would have

\begin{align*}
V_1\otimes V_3(w) \cong & \,V^f_{\omega_1+\omega_3}\oplus V^f_4 \oplus V^f_{\omega_2+\omega_7} \oplus V^f_{\omega_1+\omega_8}\oplus V^f_6\oplus V^f_{\omega_1+\omega_2}
\oplus V^f_3\oplus  V^f_{\omega_1+\omega_7}\oplus V^f_8\oplus V^f_2\oplus\\
& V^f_{2\omega_1}\oplus V^f_7\oplus V^f_1\oplus V^f_{2\omega_1+\omega_7}
\oplus V^f_{2\omega_7}\oplus V^f_{3\omega_1}\oplus \blb C\oplus \rm{mult}
\end{align*} 

From Table 2 we see that, at $w=q^4$, this tensor product has the ECC of $V_4$ and
that none of the other components with height-2 highest weight can have this ECC (we use the definition 
height$(\lambda)=\sum \lambda_i$, for a  dominant weight $\lambda = \sum \lambda_i\omega_i$).
We still need to check that this is also true for the height-3 highest weight components. 
The ECC at $V_1$ for these components are of the form  
$\xi_{1i}(z/w_1)\xi_{1j}(z/w_2)\xi_{1k}(z/w_3)$. We used the computer software 
{\it Mathematica} to check that such products (for the $i,j,k$ we need) will never coincide
with $\xi_{14}(u)$. Therefore we must have an affinization of
$V_4$ as a subrepresentation 
of $V_1\otimes V_3(q^4)$. Again assume that $V_4$ has the longest 
possible ``tail''. Then

\begin{align*}
V_1\otimes V_4(w) \cong & (V^f_1\otimes V^f_4) \oplus 
(V^f_1\otimes V^f_{\omega_2+\omega_7}) \oplus (V^f_1\otimes V^f_{\omega_1+\omega_8})\oplus 
(V^f_1\otimes V^f_6)\oplus (V^f_1\otimes V^f_{\omega_1+\omega_2})\oplus  \\ 
& (V^f_1\otimes V^f_3)\oplus (V^f_1\otimes V^f_{\omega_1+\omega_7})\oplus (V^f_1\otimes V^f_8)
\oplus (V^f_1\otimes V^f_2)\oplus  (V^f_1\otimes V^f_{2\omega_1})\oplus \\
& (V^f_1\otimes V^f_7) 
\oplus (V^f_1\otimes V^f_1) \oplus (V^f_1\otimes V^f_{2\omega_1+\omega_7})
\oplus (V^f_1\otimes V^f_{2\omega_7})\oplus (V^f_1\otimes V^f_{3\omega_1})\oplus \\
& V^f_1 \oplus 
V_4\oplus \rm{mult}
\end{align*}
Using Table 2 and the computer one checks that all terms after $V^f_1\otimes V^f_4$ never have
ECC of $V_5$. Now
$$V^f_1\otimes V^f_4 \cong V^f_{\omega_1+\omega_4}\oplus V^f_5\oplus V^f_{\omega_3+\omega_7}
\oplus V^f_{\omega_2+\omega_8}\oplus V^f_{\omega_1+\omega_6}\oplus V^f_{\omega_7+\omega_8}
\oplus V^f_{\omega_1+\omega_3}\oplus V^f_4\oplus V^f_{\omega_2+\omega_7}
\oplus V^f_{\omega_1+\omega_8}\oplus V^f_6\oplus V^f_3$$
Use Table 2 again to conclude that we obtain an affinization of $V^f_5$ as 
subrepresentation of $V_1\otimes V_4(q^5)$.

\section{The Matrices $M(q)$}

\vspace{.5cm}\noindent $\bullet\, A_n$

$$m_{ij}(q) = \frac{(q^i-q^{-i})(q^{h^{\vee}-j} - q^{-(h^{\vee}-j)})}{(q-q^{-1})(q^{h^{\vee}}-q^{-h^{\vee}})} = -q^{h^{\vee}}\frac{(q^i-q^{-i})(q^{h^{\vee}-j}-q^{-(h^{\vee}-j)})} {(q-q^{-1})(1-q^{2h^{\vee}})}$$
for $1\leq i\leq j\leq n$.

\vspace{.5cm}\noindent $\bullet\,B_n$

For $i\leq j<n${\small
\begin{align*}
m_{ij}(q) & =  \frac{(q^{2i}-q^{-2i})(q^{h^{\vee}-2j}+q^{-(h^{\vee}-2j)})(q+q^{-1})} {(q^2-q^{-2})(q^{h^{\vee}}+q^{-h^{\vee}})} = -q^{2h^{\vee}}\frac{(q^{h^{\vee}}-q^{-h^{\vee}})(q^{2i}-q^{-2i}) (q^{h^{\vee}-2j}+q^{-(h^{\vee}-2j)})} {(q-q^{-1})(1-q^{4h^{\vee}})}\\ \\
m_{in}(q) & = \frac{(q^{2i}-q^{-2i})(q+q^{-1})} {(q^2-q^{-2})(q^{h^{\vee}}+q^{-h^{\vee}})} = -q^{2h^{\vee}}\frac{(q^{h^{\vee}}-q^{-h^{\vee}})(q^{2i}-q^{-2i})} {(q-q^{-1})(1-q^{4h^{\vee}})}\\ \\
m_{nn}(q) & = \frac{(q^{h^{\vee}+1}-q^{-(h^{\vee}+1)})} {(q^2-q^{-2})(q^{h^{\vee}}+q^{-h^{\vee}})} = -q^{2h^{\vee}}\,\frac{(q^{h^{\vee}}-q^{-h^{\vee}})(\sum_{k=0}^{h^{\vee}} (-1)^kq^{h^{\vee}-2k})} {(q-q^{-1})(1-q^{4h^{\vee}})}
\end{align*}} 

\pagebreak\vspace{.5cm}\noindent $\bullet\,C_n$

$$m_{ij}(q) = \frac{(q^i-q^{-i})(q^{h^{\vee}-j} + q^{-(h^{\vee}-j)})}{(q-q^{-1})(q^{h^{\vee}}+q^{-h^{\vee}})} = -q^{2h^{\vee}}\frac{(q^{h^{\vee}}-q^{-h^{\vee}})(q^i-q^{-i})(q^{h^{\vee}-j}+q^{-(h^{\vee}-j)})} {(q-q^{-1})(1-q^{4h^{\vee}})}$$
for $1\leq i\leq j\leq n$.

\vspace{.5cm}\noindent $\bullet\,D_n$

For $i,j<n-1$,
\begin{align*}
& m_{ij}(q) = \frac{(q^i-q^{-i})(q^{n-1-j}+q^{-(n-1-j)})} {(q-q^{-1})(q^{n-1}+q^{-(n-1)})} = -q^{h^{\vee}}\,\frac{(q^{h^{\vee}/2}-q^{-h^{\vee}/2})(q^i-q^{-i}) (q^{h^{\vee}/2-j}+q^{-(h^{\vee}/2-j)})}{(q-q^{-1})(1-q^{2h^{\vee}})}\\ \\
& m_{i\,n-1}=m_{in} = \frac{(q^i-q^{-i})} {(q-q^{-1})(q^{n-1}+q^{-(n-1)})} = -q^{h^{\vee}}\,\frac{(q^{h^{\vee}/2}-q^{-h^{\vee}/2})(q^i-q^{-i})} {(q-q^{-1})(1-q^{2h^{\vee}})}\\ \\
& m_{n-1n}= \frac{(q^{n-2}-q^{-(n-2)})} {(q-q^{-1})(q+q^{-1})(q^{n-1}+q^{-(n-1)})} = -q^{h^{\vee}}\,\frac{(q^{h^{\vee}/2}-q^{-h^{\vee}/2})(q^{h^{\vee}/2-1}-q^{-(h^{\vee}/2-1)})} {(q-q^{-1})(q+q^{-1})(1-q^{2h^{\vee}})}\\ \\
& m_{n-1\,n-1}= m_{nn} = \frac{(q^{n}-q^{-n})} {(q-q^{-1})(q+q^{-1})(q^{n-1}+q^{-(n-1)})} = -q^{h^{\vee}}\,\frac{(q^{h^{\vee}/2}-q^{-h^{\vee}/2})(q^{h^{\vee}/2+1}-q^{-(h^{\vee}/2+1)})} {(q-q^{-1})(q+q^{-1})(1-q^{2h^{\vee}})}
\end{align*}                                     

\vspace{.5cm}\noindent $\bullet\,E_6$

$$\det B(q) = q^6+q^4-1+q^{-4}+q^{-6} \quad\quad\frac{1}{\det B(q)} = \frac{-q^5+2q^7-2q^9+q^{11}+q^{13}-2q^{15}+2q^{17}-q^{19}}{(q-q^{-1})(1-q^{24})}$$
The entries of $M(q)$ are of the form
$$m_{ij}(q) = \frac{n_{ij}(q)}{p(q)} \quad\quad\text{where}\quad\quad p(q) = (q-q^{-1})(1-q^{24})$$
We now list the numerators $n_{ij}$ 

{\small\begin{align*}
& n_{11}(q)=n_{55}(q) = -1+q^2-q^6+q^8+q^{16}-q^{18}+q^{22}-q^{24}\\
& n_{12}(q)=n_{21}(q) = n_{45}(q)=n_{54}(q) = -q+q^3-q^{5}+q^{9}+q^{15}-q^{19}+q^{21}-q^{23}\\
& n_{13}(q)=n_{31}(q) = n_{35}(q)=n_{53}(q) = n_{26}(q)=n_{62}(q) = n_{46}(q)=n_{64}(q) = -q^2+q^{10}+q^{14}-q^{22}\\
& n_{14}(q)=n_{41}(q) = n_{25}(q)=n_{52}(q) = -q^3+q^7-q^{9}+q^{11}+q^{13}-q^{15}+q^{17}-q^{21}\\
& n_{15}(q)=n_{51}(q) = -q^4+q^6-q^{10}+2q^{12}-q^{14}+q^{18}-q^{20}\\
& n_{16}(q)=n_{61}(q) = n_{56}(q)=n_{65}(q) = -q^3+q^5-q^{7}+q^{9}+q^{15}-q^{17}+q^{19}-q^{21}\\
& n_{22}(q)=n_{44}(q) = -1-q^6+q^{8}+q^{10}+q^{14}+q^{16}-q^{18}-q^{24}\\
& n_{23}(q)=n_{32}(q) = n_{34}(q) = n_{43}(q) = -q-q^3+q^{9}+q^{11}+q^{13}+q^{15}-q^{21}-q^{23}\\
& n_{24}(q)=n_{42}(q) = -q^2-q^4+q^{6}+2q^{12}+q^{18}-q^{20}-q^{22}\\
& n_{33}(q)=  -1-q^2-q^4+q^8+q^{10}+2q^{12}+q^{14}+q^{16}-q^{20}-q^{22}-q^{24}\\
& n_{36}(q)=n_{63}(q) = -q-q^5+q^7+q^{11}+q^{13}+q^{17}-q^{19}-q^{23}\\
& n_{66}(q)=  -1+q^2-q^4+q^8-q^{10}+2q^{12}-q^{14}+q^{16}-q^{20}+q^{22}-q^{24}
\end{align*}}

\vspace{.5cm}\noindent $\bullet\,E_7$

$$\det B(q) = q^7+q^5-q-q^{-1}+q^{-5}+q^{-7}$$
$$\frac{1}{\det B(q)} = \frac{-q^6+2q^8-2q^{10}+q^{12}+q^{24}-2q^{26}+2q^{28}+q^{34}-q^{36}}{(q-q^{-1}) (1-q^{36})}$$
The numerators of $m_{ij}(q)$ are

{\small \begin{align*}
& n_{11}(q) = -1+q^2-q^8+q^{10}-q^{16}+2q^{18}-q^{20}+q^{26}-q^{28}+q^{34}-q^{36} \\
& n_{12}(q)=n_{21}(q) = -q+q^3-q^7+q^{11}-q^{15}+q^{17}+q^{19}-q^{21}+q^{25}-q^{29}+q^{33}-q^{35}\\ 
& n_{13}(q)=n_{31}(q) = -q^2+q^4-q^6+q^{12}-q^{14}+q^{16}+q^{20}-q^{22}+q^{24}-q^{30}+q^{32}-q^{34} \\
& n_{14}(q)=n_{41}(q) = n_{27}(q)=n_{72}(q) = n_{36}(q)=n_{63}(q) = -q^3+q^{15}+q^{21}-q^{33}\\ 
& n_{15}(q)=n_{51}(q) = n_{26}(q)=n_{26}(q) = -q^4+q^8-q^{10}+q^{14}+q^{22}-q^{26}+q^{28}-q^{32}\\ 
& n_{16}(q)=n_{61}(q) = -q^5+q^7-q^{11}+q^{13}+q^{23}-q^{25}+q^{29}-q^{31}\\ 
& n_{17}(q)=n_{71}(q) = -q^4+q^6-q^{8}+q^{10}-q^{12}+q^{14}+q^{22}-q^{24}+q^{26}-q^{28}+q^{30}-q^{32} \\
& n_{22}(q) = -1+q^4-q^6-q^{8}+q^{10}+q^{12}-q^{14}+2q^{18}-q^{22}+q^{24}+q^{26}-q^{28}-q^{30}+q^{32}-q^{36}\\
& n_{23}(q)=n_{32}(q) = -q-q^7+q^{11}+q^{17}+q^{19}+q^{25}-q^{29}-q^{35}\\ 
& n_{24}(q)=n_{42}(q) = n_{35}(q)=n_{53}(q) = -q^2-q^4+q^{14}+q^{16}+q^{20}+q^{22}-q^{32}-q^{34}\\ 
& n_{25}(q)=n_{52}(q) = -q^3-q^5+q^7-q^{11}+q^{13}+q^{15}+q^{21}+q^{23}-q^{25}+q^{29}-q^{31}-q^{33}\\
& n_{33}(q) = -1-q^4-q^8+q^{10}+q^{14}+2q^{18}+q^{22}+q^{26}-q^{28}-q^{32}-q^{36}\\
& n_{34}(q)=n_{43}(q) = -q-q^3-q^5+q^{13}+q^{15}+q^{17}+q^{19}+q^{21}+q^{23}-q^{31}-q^{33}-q^{35}\\
& n_{37}(q)=n_{73}(q) = -q^2-q^6+q^{8}-q^{10}+q^{12}+q^{16}+q^{20}+q^{24}-q^{26}+q^{28}-q^{30}-q^{34}\\
\end{align*}
\begin{align*}
& n_{44}(q) = -1-q^2-q^4-q^6+q^{12}+q^{14}+q^{16}+2q^{18}+q^{20}+q^{22}+q^{24}-q^{30}-q^{32}-q^{34}-q^{36}\\
& n_{45}(q)=n_{54}(q) = -q-q^3-q^7+q^{11}+q^{15}+q^{17}+q^{19}+q^{21}+q^{25}-q^{29}-q^{33}-q^{35}\\
& n_{46}(q)=n_{64}(q) = n_{57}(q)=n_{75}(q) = -q^2-q^8+q^{10}+q^{16}+q^{20}+q^{26}-q^{28}-q^{34}\\
& n_{47}(q)=n_{74}(q) = -q-q^5+q^{13}+q^{17}+q^{19}+q^{23}-q^{31}-q^{35}\\
& n_{55}(q) = -1-q^6+q^{12}+2q^{18}+q^{24}-q^{30}-q^{36}\\
& n_{56}(q)=n_{65}(q) = -q+q^3-q^5+q^{13}-q^{15}+q^{17}+q^{19}-q^{21}+q^{23}-q^{31}+q^{33}-q^{35}\\
& n_{66}(q) = -1+q^2-q^6+q^{8}-q^{10}+q^{12}-q^{16}+2q^{18}-q^{20}+q^{24}-q^{26}+q^{28}-q^{30}+q^{34}-q^{36}\\
& n_{67}(q)=n_{76}(q) = -q^3+q^5-q^7+q^{11}-q^{13}+q^{15}+q^{21}-q^{23}+q^{25}-q^{29}+q^{31}-q^{33}\\
& n_{77}(q) = -1+q^2-q^4+q^{14}-q^{16}+2q^{18}-q^{20}+q^{22}-q^{32}+q^{34}-q^{36}
\end{align*}}

\vspace{.5cm}\noindent $\bullet\,E_8$

$$\det B(q) = q^8+q^6-q^2-1-q^{-2}+q^{-6}+q^{-8} \qquad\qquad \frac{1}{\det B(q)} = \frac{p(q)}{(q-q^{-1})(1-q^{60})}$$ where 
\begin{align*}
p(q) = & -q^7+2q^9-2q^11+q^{13}-q^{17}+2q^{19}-2q^{21}+q^{23}\\
&+q^{37}-2q^{39}+2q^{41}-2q^{43} +q^{47}-2q^{49}+2q^{51}-q^{53}
\end{align*}
The numerators of $m_{ij}(q)$ are 

{\small
\begin{align*}
& n_{11} = -1+q^2-q^{10}+q^{12}-q^{18}+q^{20}-q^{28}+ 2q^{30}-q^{32}+q^{40}-q^{42}+q^{48}-q^{50} +q^{58}-q^{60}\\
& n_{12} = -q+q^3-q^{9}+q^{13}-q^{17}+q^{21}-q^{27}+q^{29} +q^{31}-q^{33}+q^{39}-q^{43}+q^{47}-q^{51} +q^{57}-q^{59}\\
& n_{13} = -q^2+q^4-q^{8}+q^{14}-q^{16}+q^{22}-q^{26}+q^{28} +q^{32}-q^{34}+q^{38}-q^{44}+q^{46}-q^{52}+q^{56}-q^{58}\\
& n_{14}=n_{78} = -q^3+q^{5}-q^{7}+q^{23}-q^{25}+q^{27} +q^{33}-q^{35}+q^{37}-q^{53}+q^{55}-q^{57}\\
& n_{15}=n_{28}=n_{37} = -q^4-q^{14}+q^{16}+q^{26} +q^{34}+q^{44}-q^{46}-q^{56}\\
& n_{16}=n_{27} = -q^5+q^{9}-q^{11}+q^{19}-q^{21}+q^{25} +q^{35}-q^{39}+q^{41}-q^{49}+q^{51}-q^{55}\\
& n_{17} = -q^6+q^{8}-q^{12}+q^{14}-q^{16}+q^{18}-q^{22}+q^{24} +q^{36}-q^{38}+q^{42}-q^{44}+q^{46}-q^{48}+q^{52}-q^{54}\\
& n_{18}= -q^5+q^7-q^{9}+q^{11}-q^{13}+q^{17}-q^{19}+q^{21}-q^{23}+q^{25} \\
& \hspace{1cm} +q^{35}-q^{37}+q^{39}-q^{41}+q^{43}-q^{47}+q^{49}-q^{51}+q^{53}-q^{55}\\
& n_{22} = -1+q^4-q^8-q^{10}+q^{12}+q^{14}-q^{16}-q^{18}+q^{20}+q^{22}-q^{26}+ 2q^{30}\\
& \hspace{1cm} -q^{34}+q^{38}+q^{40}-q^{42}-q^{44}+q^{46}+q^{48}-q^{50}-q^{52}+q^{56}-q^{60}\\
& n_{23} = -q+q^5-q^{7}-q^9+q^{13}-q^{17}+q^{21}+q^{23}-q^{25}+q^{29}\\
& \hspace{1cm} +q^{31}-q^{35}+q^{37}+q^{39}-q^{43}+q^{47}-q^{51}-q^{53}+q^{55}-q^{59}\\
& n_{24}=n_{57}=n_{68} = -q^2-q^{8}+q^{22}+q^{28} +q^{32}+q^{38}-q^{52}-q^{58}\\
& n_{25}=n_{36} = -q^3-q^{5}-q^{13}+q^{17}+q^{25}+q^{27} +q^{33}+q^{35}+q^{43}-q^{47}-q^{55}-q^{57}\\
& n_{26} = -q^4-q^6+q^{8}-q^{12}+q^{18}-q^{22}+q^{24}+q^{26} +q^{34}+q^{36}-q^{38}+q^{42}-q^{48}+q^{52}-q^{54}-q^{56}\\
\end{align*}
\begin{align*}
& n_{33} = -1-q^8-q^{10}+q^{12}-q^{18}+q^{22}+2q^{30}
+q^{38}+q^{40}-q^{42}+q^{48}-q^{50}-q^{52}-q^{60}\\
& n_{34}=n_{58} = -q-q^5-q^{9}+q^{21}+q^{25}+q^{29} +q^{31}+q^{35}+q^{39}-q^{51}-q^{55}-q^{59}\\
& n_{35} = -q^2-q^4-q^{6}-q^{12}+q^{18}+q^{24}+q^{26}+q^{28}
+q^{32}+q^{34}+q^{36}+q^{42}-q^{48}-q^{54}-q^{56}-q^{58}\\
& n_{38}= -q^3-q^7+q^{9}-q^{11}+q^{19}-q^{21}+q^{23}+q^{27} +q^{33}+q^{37}-q^{39}+q^{41}-q^{49}+q^{51}-q^{53}-q^{57}\\
& n_{44} = -1-q^4-q^{6}-q^{10}+q^{20}+q^{24}+q^{26}+2q^{30}
+q^{34}+q^{36}+q^{40}-q^{50}-q^{54}-q^{56}-q^{60}\\
& n_{45}= -q-q^3-q^5-q^{7}-q^{11}+q^{19}+q^{23}+q^{25}+q^{27}+q^{29}\\
& \hspace{1cm} +q^{31}+q^{33}+q^{35}+q^{37}+q^{41}-q^{49}-q^{53}-q^{55}-q^{57}-q^{59}\\
& n_{46} = -q^2-q^4-q^{8}-q^{14}+q^{16}+q^{22}+q^{26}+q^{28}
+q^{32}+q^{34}+q^{38}+q^{44}-q^{46}-q^{52}-q^{56}-q^{58}\\
& n_{47}= -q^3-q^9+q^{11}-q^{13}+q^{17}-q^{19}+q^{21}+q^{27} +q^{33}+q^{39}-q^{41}+q^{43}-q^{47}+q^{49}-q^{51}-q^{57}\\
& n_{48} = -q^2-q^6-q^{12}+q^{14}-q^{16}+q^{18}+q^{24}+q^{28}
+q^{32}+q^{36}+q^{42}-q^{44}+q^{46}-q^{48}-q^{54}-q^{58}\\
& n_{55} = -1-q^2-q^4-q^{6}-q^8-q^{10}+q^{20}+q^{22}+q^{24}+q^{26}+q^{28}+2q^{30}\\
& \hspace{1cm} +q^{32}+q^{34}+q^{36}+q^{38}+q^{40}-q^{50}-q^{52}-q^{54}-q^{56}-q^{58}-q^{60}\\
& n_{56}= -q-q^3-q^7-q^{9}+q^{21}+q^{23}+q^{27}+q^{29} +q^{31}+q^{33}+q^{37}+q^{39}-q^{51}-q^{53}-q^{57}-q^{59}\\
& n_{66} = -1-q^6-q^{10}+q^{14}+q^{16}+q^{20}+q^{24}+2q^{30}
+q^{36}+q^{40}-q^{44}+q^{46}-q^{50}-q^{54}-q^{60}\\
& n_{67}= -q+q^3-q^5-q^{11}+q^{13}-q^{17}+q^{19}+q^{25}-q^{27}+q^{29}\\
& \hspace{1cm} +q^{31}-q^{33}+q^{35}+q^{41}-q^{43}+q^{47}-q^{49}-q^{55}+q^{57}-q^{59}\\
& n_{77} = -1+q^2-q^6+q^{8}-q^{10}+q^{14}-q^{16}+q^{20}-q^{22}+q^{24}-q^{28}+2q^{30}\\
& \hspace{1cm} -q^{32}+q^{36}-q^{38}+q^{40}-q^{44}+q^{46}-q^{50}+q^{52}-q^{54}+q^{58}-q^{60}\\
& n_{88} = -1+q^2-q^4-q^{10}+q^{12}-q^{14}+q^{16}-q^{18}+q^{20}+q^{26}-q^{28}+2q^{30}\\
& \hspace{1cm} -q^{32}+q^{34}+q^{40}-q^{42}+q^{44}-q^{46}+q^{48}-q^{50}-q^{56}+q^{58}-q^{60}\\
\end{align*}

\vspace{.5cm}\noindent $\bullet\,F_4$

$$B(q) = 
\begin{pmatrix}
[2]_q & -1 & 0 & 0\\
-1 & [2]_q & -[2]_q & 0\\
0 & -[2]_q & [4]_q & -[2]_q\\
0 & 0 & -[2]_q & [4]_q
\end{pmatrix}
\quad\quad\det B(q) = q^8+2q^6+q^4-q^2-2-q^{-2}+q^{-4}+2q^{-6}+q^{-8}$$ 
$$\frac{1}{\det B(q)} = \frac{-q^6+2q^8-2q^{10}+q^{12}+q^{24}-2q^{26}+2q^{28}-q^{30}}{(q+q^{-1})(q-q^{-1})(1-q^{36})}$$
and the entries of $M(q)$ are of the for $m_{ij}(q) = \frac{n_{ij}(q)}{(q-q^{-1})(1-q^{36})}$, where the numerators are {\small
\begin{align*}
& n_{11}(q) = -1+q^2-q^6+q^8-q^{10}+q^{12}-q^{16}+2q^{18} -q^{20}+q^{24}-q^{26}+q^{28}-q^{30}+q^{34}-q^{36}\\
& n_{12}(q)=n_{21}(q) = -q+q^3-q^5+q^{13}-q^{15}+q^{17}+q^{19}-q^{21}+q^{23}-q^{31}+q^{33}-q^{35}\\
& n_{13}(q)=n_{31}(q) = -q^2-q^8+q^{10}+q^{16}+q^{20}+q^{26}-q^{28}-q^{34}\\
& n_{14}(q)=n_{41}(q) = -q^4+q^8-q^{10}+q^{14}+q^{22}-q^{26}+q^{28}-q^{32}\\
& n_{22}(q) = -1-q^6+q^{12}+2q^{18}+q^{24}-q^{30}-q^{36}\\
& n_{23}(q)=n_{32}(q) = -q-q^3-q^7+q^{11}+q^{15}+q^{17}+q^{19}+q^{21}+q^{25}-q^{29}-q^{33}-q^{35}\\
& n_{24}(q)=n_{42}(q) = -q^{3}-q^5+q^7-q^{11}+q^{13}+q^{15}+q^{21}+q^{23}-q^{25}+q^{29}-q^{31}-q^{33}\\
& n_{33}(q) = -1-q^2-q^4-q^6+q^{12}+q^{14}+q^{16}+2q^{18}+q^{20}+q^{22}+q^{24}-q^{30}-q^{32}-q^{34}-q^{36}\\
& n_{34}(q)=n_{43}(q) = -q^2-q^4+q^{14}+q^{16}+q^{20}+q^{22}-q^{32}-q^{34}\\
& n_{44}(q) = -1+q^4-q^6-q^8+q^{10}+q^{12}-q^{14}+2q^{18}-q^{22}+q^{24}+q^{26}-q^{28}-q^{30}+q^{32}-q^{36}
\end{align*} } 

\vspace{.5cm}\noindent $\bullet\,G_2$

$$B(q) = \begin{pmatrix}[2]_q & -[3]_q\\ -[3]_q & [6]_q\end{pmatrix}\quad\quad\quad D(q) = \begin{pmatrix} 1 & 0\\ 0 & [3]_q\end{pmatrix} \quad\quad\det B(q) = q^6+q^4-1+q^{-4}+q^{-6} $$
$$\frac{1}{\det B(q)} = -\frac {(q-q^{-1})(q^7+q^{11}+q^{13}+q^{17})}{(1-q^{24})} = \frac{q^6-q^8+q^{10}-q^{14}+q^{16}-q^{18}}{(1-q^{24})}$$
$$
M(q) = \frac{1}{\det B(q)}\,\begin{pmatrix} [6]_q & [3]_q^2\\ \\ \,[3]_q^2 & [2]_q[3]_q^2\end{pmatrix} = 
\frac{1}{\det B(q)}\,\begin{pmatrix} [6]_q & [3]_q^2\\  \\ \,[3]_q^2 & (q+q^{-1})[3]_q^2\end{pmatrix}
$$  



\end{document}